\documentclass[12pt]{article}
\usepackage{amssymb,amsmath,amsthm}
\newtheorem{theorem}{Theorem}
\newtheorem{lemma}{Lemma}
\newtheorem{proposition}{Proposition}

\newcommand{\tr}{\,\mathrm{tr}\,}
\newcommand{\ad}{\,\mathrm{ad}\,}
\begin{document}

\begin{center}

{\Large {\bf Elementary equivalence \\

\bigskip

of Chevalley groups over fields}}

\bigskip
\bigskip

{\large \bf E.~I.~Bunina}

\end{center}
\bigskip

\section*{Introduction}\leavevmode

Two models ${\cal U}$ and ${\cal U}'$ of the same first order
language~$\cal L$
(for example, two groups or two rings) are called
\emph{elementarily equivalent} if every sentence~$\varphi$ of the 
language~$\cal L$ holds in~${\cal U}$ if and only if it holds
in~${\cal U}'$. Any two finite models of the same language
are elementarily equivalent if and only if they are isomorphic.
Any two isomorphic models are elementarily equivalent, but for infinite
models the converse is not true.
For example, the field $\mathbb C$
of complex numbers and the field $\overline {\mathbb Q}$ of algebraic numbers
are elementarily equivalent, but not isomorphic, since they have different
powers (for more detailed examples see~\cite{KeisChang}).

The first results on connection between elementary properties of some 
models and elementary properties of derivative models were obtained
by A.I.~Maltsev 1961 in~\cite{Maltsev}.
He proved that the groups $G_n(K)$ and $G_m(L)$
(where $G=GL,SL,PGL,PSL$, $n,m\ge 3$, $K$,~$L$ are fields of
characteristics~$0$)
are elementarily equivalent if and only if
 $m=n$ and the fields $K$~and~$L$ are elementarily equivalent.

Studying of these problems was continued in 1992, when with 
the help of ultrapower construction and the Isomorphism theorem~\cite{KeisChang} 
K.I.~Beidar and A.V.~Mikhalev~\cite{BeiMikh} formulated the general approach
to the problems of elementary equivalence of different
algebraic structures and generalized the Maltsev theorem
to the case when $K$~and~$L$ are skewfields and associative rings.

In 1998--2004 E.I.~Bunina continued to study some problems of
this type (see~\cite{Bun1}--\cite{Bun3}). The results of
A.I.~Maltsev were generalized to
unitary linear groups over skewfields and associative rings
with involutions,
and also Chevalley groups over algebraically closed fields.

In the paper~\cite{BunChev} we announced the following results
on elementary equivalence of Chevalley groups over fields:

\smallskip

{\bf Theorem 1.} \emph{Let $G_\pi (\Phi,K)$ and $G_{\pi'}(\Phi',K')$
\emph{(}or $E_\pi (\Phi,K)$ and $E_{\pi'}(\Phi',K'))$ be two \emph{(}elementary\emph{)}
Chevalley groups over infinite fields $K$ and~$K'$ of characteristics $\ne 2$.
Then elementary equivalence of these Chevalley groups
implies $\Phi\cong \Phi'$ and $R\equiv R'$.}

\smallskip

{\bf Theorem 2.} \emph{Let $G=G_\pi (\Phi,K)$ and $G'=G_{\pi'}(\Phi,K')$
\emph{(}or $E_\pi (\Phi,K)$ and $E_{\pi'}(\Phi,K'))$ be two 
\emph{(}elementary\emph{)}
Chevalley groups over elementary equivalent fields $K$ and~$K'$,
where representations $\pi$ and $\pi'$ have the same weight lattices.
Then the groups $G$ and $G'$ are elementarily equivalent.}

\smallskip

The given work is devoted to the detailed proof
of these two theorems, and even more strict \emph{main theorem}:

\smallskip

{\bf Main theorem.} \emph{Let
 $G=G_\pi (\Phi,K)$ and $G'=G_{\pi'}(\Phi',K')$
\emph{(}or $E_\pi (\Phi,K)$ and $E_{\pi'}(\Phi',K'))$ be two
 \emph{(}elementary\emph{)}
Chevalley groups over infinite fields $K$ and~$K'$ of characteristic
$\ne 2$, with weight lattices $\Lambda$ and $\Lambda'$, respectively.
Then the groups $G$ and $G'$ are elementarily equivalent if and only if
the root systems $\Phi$ and $\Phi'$ are isomorphic, the fields $K$  and
 $K'$ are elementarily equivalent, the lattices  $\Lambda$ and $\Lambda'$ 
coincide.}

\section{Basic facts about Chevalley groups.}\leavevmode

\subsection{Root systems.}\leavevmode

More detailes about root systems and their properties can be found
in the books \cite{Hamfris},~\cite{Burbaki}.

A finite nonempty set $\Phi\subset \mathbb R^l$ of vectors
of the Euclidean space~$\mathbb R^l$ is called a \emph{root system},
if it generates $\mathbb R^l$, does not contain~$0$
and satisfies the following properties:

1) $\forall \alpha\in \Phi\ (c\cdot \alpha  \in \Phi\Leftrightarrow
c=\pm1)$;

2) if we introduce 
$$
\langle \alpha,\beta\rangle:= \frac{2(\alpha,\beta)}{(\alpha,\alpha)}
\quad \text{(\emph{the reflection coefficient})}
$$
for $\alpha,\beta\in \mathbb R^l$, then for any $\alpha,\beta\in \Phi$ we
have $\langle \alpha,\beta\rangle \in \mathbb Z$;

3) let for $\alpha\in \mathbb R^l$ $w_\alpha$ be the reflection
under a hyperplane, ortoghonal to the vector~$\alpha$, i.\,e. $\forall \beta
\in \mathbb R^l$ 
$$
w_\alpha (\beta)=\beta-\langle \alpha,\beta\rangle \alpha.
$$
Then for any $\alpha,\beta\in \Phi$ we have $w_\alpha(\beta)\in \Phi$, 
i.\,e. the set $\Phi$ is invariant
under the action of all reflections~$w_\alpha$, $\alpha\in  \Phi$.

If $\Phi$ is a root system, then its elements are called \emph{roots}.

The group $W$, generated by all reflections $w_\alpha$, 
$\alpha\in \Phi$, is called the \emph{Weil group}
of the root system~$\Phi$.

If we draw a hyperplane in the space $\mathbb R^l$ that does not contain any roots
from~$\Phi$, then all roots are divided into
two disjoint sets of \emph{positive} ($\Phi^+$) and
\emph{negative} ($\Phi^-$) roots.
The  \emph{system
of simple roots} is a set $\Delta=
\{ \alpha_1,\dots,\alpha_l\}\subset 
\Phi^+$ such that
any positive root $\beta\in \Phi^+$ can be uniquely represented in the form
$n_1\alpha_1+\dots+n_l\alpha_l$, where
$n_1,\dots,n_l\in \mathbb Z^+$.

For any root system $\Phi$ there exists a system of simple roots.
The number $l$ is called a \emph{rank} of the root system~$\Phi$.

We are mostly interested in \emph{undecomposible} root systems,
i.\,e., such systems~$\Phi$ that can not be represented 
as the disjoint union  $\Phi=\Phi_1\cup \Phi_2$ of two
sets with mutually orthogonal roots.

By a root system we can construct the following \emph{Dynkin scheme}.
It is a graph, that is constructed as follows:
its vertices are corresponded to the simple roots $\alpha_1,\dots,\alpha_l$, 
the vertices with the numbers 
$i$ and~$j$ are connected by the edge if
 $\langle \alpha_i,\alpha_l\rangle\ne 0$.
If $|\alpha_i|=|\alpha_j|$, then $\langle \alpha_i,\alpha_j\rangle=
\langle \alpha_j,\alpha_i\rangle$ and the number of edges between the vertices
 $i$ and~$j$ is $|\langle \alpha_i,\alpha_j\rangle|$.  If $|\alpha_i|>
|\alpha_j|$, then $\langle \alpha_i,\alpha_j\rangle< \langle \alpha_j,\alpha_i\rangle$
and $|\langle \alpha_i,\alpha_j\rangle|=1$. In this case there are 
 $|\langle \alpha_j,\alpha_i\rangle|$ edges between $i$ and~$j$
and we put an arrow from the long root to the short one.

The Dynkin scheme and the root system are uniquely corresponded to each other.

All indecomposible root systems up to isomorphism are divided to $4$
infinite (\emph{classical}) 
series $A_l$ ($l\ge 1$), $B_l$ ($l\ge 2$), $C_l$ ($l\ge 3$) and
$D_l$ ($l\ge 4$) and $5$ separate (\emph{exceptional}) 
cases $E_6$, $E_7$, $E_8$, $F_4$,
and~$G_2$.

These are Dynkin schemes for the root systems:

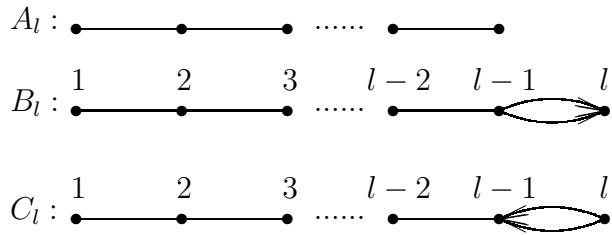
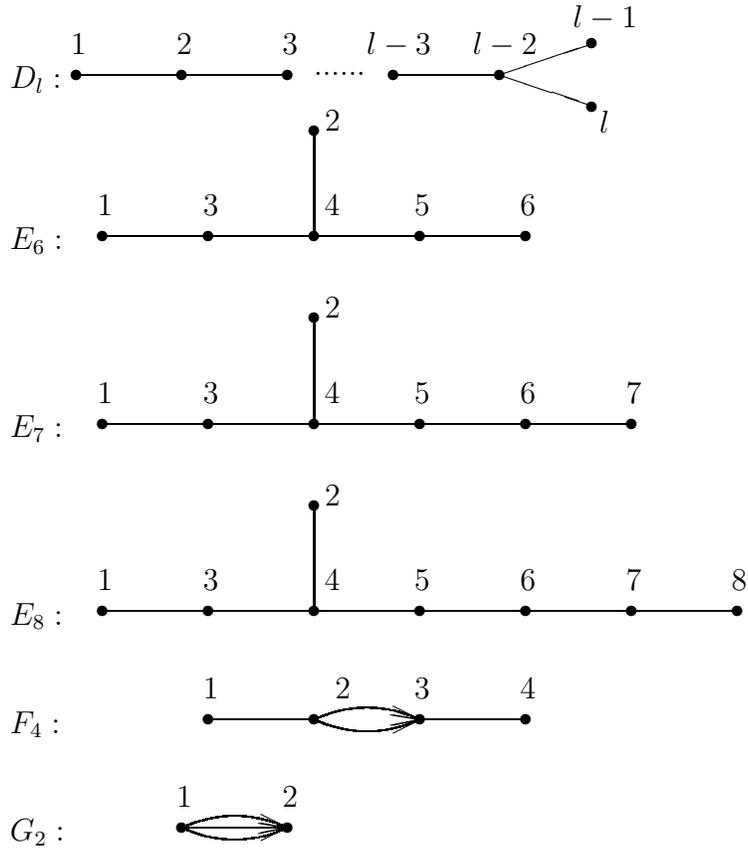
\begin{picture}(200,35)
\put(5,20){$A_l:$}
\put(30,20){\line(1,0){40}}
\put(30,20){\circle*{4}}
\put(70,20){\line(1,0){40}}
\put(70,20){\circle*{4}}
\put(110,20){\circle*{4}}
\put(150,20){\circle*{4}}
\put(190,20){\circle*{4}}
\put(120,20){......}
\put(150,20){\line(1,0){40}}
\end{picture}

\begin{picture}(250,30)
\put(5,20){$B_l:$}
\put(30,20){\line(1,0){40}}
\put(30,20){\circle*{4}}
\put(70,20){\line(1,0){40}}
\put(70,20){\circle*{4}}
\put(110,20){\circle*{4}}
\put(150,20){\circle*{4}}
\put(190,20){\circle*{4}}
\put(230,20){\circle*{4}}
\put(120,20){......}
\put(150,20){\line(1,0){40}}
\qbezier(190,20)(210,29)(230,20)
\qbezier(190,20)(210,11)(230,20)
\put(230,20){\line(-2,1){10}}
\put(230,20){\line(-2,-1){10}}
\put(230,20){\line(-6,1){11}}
\put(230,20){\line(-5,-1){11}}
\put(28,28){$1$}
\put(68,28){$2$}
\put(108,28){$3$}
\put(140,28){$l-2$}
\put(180,28){$l-1$}
\put(228,28){$l$}
\end{picture}

\begin{picture}(250,40)
\put(5,20){$C_l:$}
\put(30,20){\line(1,0){40}}
\put(30,20){\circle*{4}}
\put(70,20){\line(1,0){40}}
\put(70,20){\circle*{4}}
\put(110,20){\circle*{4}}
\put(150,20){\circle*{4}}
\put(190,20){\circle*{4}}
\put(230,20){\circle*{4}}
\put(120,20){......}
\put(150,20){\line(1,0){40}}
\qbezier(190,20)(210,29)(230,20)
\qbezier(190,20)(210,11)(230,20)
\put(190,20){\line(2,1){10}}
\put(190,20){\line(2,-1){10}}
\put(190,20){\line(6,1){11}}
\put(190,20){\line(5,-1){11}}
\put(28,28){$1$}
\put(68,28){$2$}
\put(108,28){$3$}
\put(140,28){$l-2$}
\put(180,28){$l-1$}
\put(228,28){$l$}
\end{picture}

\begin{picture}(250,40)
\put(5,25){$D_l:$}
\put(30,30){\line(1,0){40}}
\put(30,30){\circle*{4}}
\put(70,30){\line(1,0){40}}
\put(70,30){\circle*{4}}
\put(110,30){\circle*{4}}
\put(150,30){\circle*{4}}
\put(190,30){\circle*{4}}
\put(225,18){\circle*{4}}
\put(225,42){\circle*{4}}
\put(120,30){......}
\put(150,30){\line(1,0){40}}
\put(190,30){\line(3,1){35}}
\put(190,30){\line(3,-1){35}}
\put(28,38){$1$}
\put(68,38){$2$}
\put(108,38){$3$}
\put(140,38){$l-3$}
\put(180,38){$l-2$}
\put(218,48){$l-1$}
\put(228,8){$l$}
\end{picture}

\begin{picture}(250,60)
\put(5,25){$E_6:$}
\put(40,30){\line(1,0){40}}
\put(40,30){\circle*{4}}
\put(80,30){\line(1,0){40}}
\put(80,30){\circle*{4}}
\put(120,30){\circle*{4}}
\put(120,70){\circle*{4}}
\put(160,30){\circle*{4}}
\put(200,30){\circle*{4}}
\put(120,30){\line(1,0){40}}
\put(120,30){\line(0,1){40}}
\put(160,30){\line(1,0){40}}
\put(38,38){$1$}
\put(78,38){$3$}
\put(124,38){$4$}
\put(124,70){$2$}
\put(158,38){$5$}
\put(198,38){$6$}
\end{picture}

\begin{picture}(250,70)
\put(5,25){$E_7:$}
\put(40,30){\line(1,0){40}}
\put(40,30){\circle*{4}}
\put(80,30){\line(1,0){40}}
\put(80,30){\circle*{4}}
\put(120,30){\circle*{4}}
\put(120,70){\circle*{4}}
\put(160,30){\circle*{4}}
\put(200,30){\circle*{4}}
\put(240,30){\circle*{4}}
\put(120,30){\line(1,0){40}}
\put(120,30){\line(0,1){40}}
\put(160,30){\line(1,0){40}}
\put(200,30){\line(1,0){40}}
\put(38,38){$1$}
\put(78,38){$3$}
\put(124,38){$4$}
\put(124,70){$2$}
\put(158,38){$5$}
\put(198,38){$6$}
\put(238,38){$7$}
\end{picture}

\begin{picture}(250,70)
\put(5,25){$E_8:$}
\put(40,30){\line(1,0){40}}
\put(40,30){\circle*{4}}
\put(80,30){\line(1,0){40}}
\put(80,30){\circle*{4}}
\put(120,30){\circle*{4}}
\put(120,70){\circle*{4}}
\put(160,30){\circle*{4}}
\put(200,30){\circle*{4}}
\put(240,30){\circle*{4}}
\put(280,30){\circle*{4}}
\put(120,30){\line(1,0){40}}
\put(120,30){\line(0,1){40}}
\put(160,30){\line(1,0){40}}
\put(200,30){\line(1,0){40}}
\put(240,30){\line(1,0){40}}
\put(38,38){$1$}
\put(78,38){$3$}
\put(124,38){$4$}
\put(124,70){$2$}
\put(158,38){$5$}
\put(198,38){$6$}
\put(238,38){$7$}
\put(278,38){$8$}
\end{picture}

\begin{picture}(250,40)
\put(5,25){$F_4:$}
\put(80,30){\line(1,0){40}}
\put(80,30){\circle*{4}}
\put(120,30){\circle*{4}}
\put(160,30){\circle*{4}}
\put(200,30){\circle*{4}}
\put(160,30){\line(1,0){40}}
\put(78,38){$1$}
\put(128,38){$2$}
\put(158,38){$3$}
\put(198,38){$4$}
\qbezier(120,30)(140,39)(160,30)
\qbezier(120,30)(140,21)(160,30)
\put(160,30){\line(-2,1){10}}
\put(160,30){\line(-2,-1){10}}
\put(160,30){\line(-6,1){11}}
\put(160,30){\line(-5,-1){11}}

\end{picture}

\begin{picture}(250,40)
\put(5,25){$G_2:$}
\put(70,30){\circle*{4}}
\put(110,30){\circle*{4}}
\put(70,30){\line(1,0){40}}
\put(68,38){$1$}
\put(108,38){$2$}
\qbezier(70,30)(90,39)(110,30)
\qbezier(70,30)(90,21)(110,30)
\put(110,30){\line(-2,1){10}}
\put(110,30){\line(-2,-1){10}}
\put(110,30){\line(-6,1){11}}
\put(110,30){\line(-5,-1){11}}
\end{picture}

\subsection{Semisimple Lie algebras.}\leavevmode

More details about semisimple Lie algebras
can be found in the book~\cite{Hamfris}.

\emph{Lie algebra} $\cal L$ over a field~$K$ 
is a linear space over~$K$,
with a binary operation 
$x,y\mapsto [x,y]$ (\emph{Lie bracket}),
that is linear by the both variables and satisfies the following conditions:

1) \emph{anti-commutativity}:
$\forall x,y\in {\cal L}\ [x,y]=-[y,x];$

2) \emph{Jacobi identity}:
$\forall x,y,z\in {\cal L}\ [x,[y,z]]+[y,[z,x]]+[z,[x,y]]=0.$

The dimension of Lie algebra is defined 
by its dimension as a linear space over~$K$.
So a Lie algebra is called \emph{finitely dimensional},
if the space $\cal L$ is finitely dimensional.

A subspace ${\cal L}'$ of~$\cal L$ 
as a linear space is called a
\emph{subalgebra} of~$\cal L$ if
$\forall x,y\in {\cal L}'$ $[x,y]\in {\cal L}'$.
A subalgebra ${\cal L}'$ of~$\cal L$ is called its \emph{ideal} if
$\forall x\in {\cal L}'$ $\forall y\in {\cal L}$\ $[x,y]\in {\cal L}'$.

The \emph{commutant} of Lie algebras ${\cal L}_1$ and ${\cal L}_2$, 
that are subalgebras of the same Lie algebra~$\cal L$, 
is its subalgebra
${\cal L}'=[{\cal L}_1,{\cal L}_2],$
generated by all elements $[x,y]$ for $x\in {\cal L}_1$, $y\in {\cal L}_2$.

The sequence
$$
{\cal L}^0={\cal L}, {\cal L}^1=[{\cal L},{\cal L}^0], {\cal L}^2=[
{\cal L},{\cal L}^1],\dots, {\cal L}^{n+1}=[{\cal L},{\cal L}^n],\dots
$$
is called the \emph{lower central series} of~${\cal L}$.
If for some $n\in \mathbb N$ we have ${\cal L}^n=0$, then
$\cal L$ is called \emph{nilpotent}. 

The sequence
$$
{\cal L}^{(0)}={\cal L}, {\cal L}^{(1)}=[{\cal L}^{(0)},{\cal L}^{(0)}], 
{\cal L}^{(2)}=[{\cal L}^{(1)},{\cal L}^{(1)}],\dots, 
{\cal L}^{(n+1)}=[{\cal L}^{(n)},{\cal L}^{(n)}],\dots
$$
is called the \emph{derivative series} of~$\cal L$. If for some
$n\in \mathbb N$ we have 
${\cal L}^{(n)}=0$, then~$\cal L$ is called \emph{solvable}.

Let us consider a finitely dimensional Lie algebra~$\cal L$.
The greatest solvable ideal of~$\cal L$,
that is the sum of all its solvable ideals,
is called the \emph{radical} of~$\cal L$. A Lie algebra with
zero radical is called \emph{semisimple}.
Noncommutative Lie algebra  $\cal L$ is called  \emph{simple}, 
if it contains exactly two ideals:
$0$ and $\cal L$. 

A finitely dimensional semisimple Lie algebra over $\mathbb C$
is a direct sum of simple Lie algebras.

The \emph{normalizer} of a subalgebra ${\cal L}'$ in a Lie algebra~$\cal L$ 
is a subalgebra
$$
N_{\cal L}({\cal L}'):=\{ x\in {\cal L}\,|\, \forall y\in {\cal L}' \ [x,y]\in 
{\cal L}'\}.
$$

A \emph{Cartan subalgebra} of a Lie algebra~${\cal L}$ is 
a nilpotent self-normalized subalgebra~$\cal H$. For a semisimple
Lie algebra
it is Abelian and is defined up to an automorphism of~$\cal L$.

Let $\cal L$ be semisimple finitely dimensional Lie algebra
over~$\mathbb C$, $\cal H$ be its Cartan subalgebra.
Consider the space~${\cal H}^*$.
Let for $\alpha\in {\cal H}^*$
$$
{\cal L}_\alpha:=\{ x\in {\cal L}\mid [h,x]=\alpha(h)x\text{ for any }
h\in {\cal H}\}.
$$

In this case ${\cal L}_0={\cal H}$, and the algebra $\cal L$ 
can be decomposed as ${\cal L}={\cal H}
\oplus \sum\limits_{\alpha\ne 0} 
{\cal L}_\alpha$, and if ${\cal L}_\alpha\ne 0$, then $\dim {\cal L}_\alpha=1$,
all such nonzero $\alpha\in {\cal H}$ that ${\cal L}_\alpha\ne 0$,
form some root system~$\Phi$.
A root system $\Phi$ and a semisimple Lie algebra $\cal L$ over~$\mathbb C$ 
uniquely define each other.

On a Lie algebra $\cal L$ we can introduce the bilinear \emph{Killing form}
$$
\varkappa(x,y)=\tr (\ad x\ad y),
$$
where $\ad x\in GL({\cal L})$, $\ad x: z\mapsto [x,z]$ is the \emph{adjoint
representation}. For a semisimple Lie algebra the restriction of the
Killing form to~$\cal H$ is not degenerate,
so we can identify the spaces $\cal H$ and ${\cal H}^*$.

We can choose a basis $\{ h_1, \dots, h_l\}$ of~$\cal H$ and for every
$\alpha\in \Phi$ elements $x_\alpha \in {\cal L}_\alpha$ so that

1) $\{ h_i; x_\alpha\}$ form a basis in~$\cal L$;

2) $[h_i,h_j]=0$;

3) $[h_i,x_\alpha]=\langle \alpha,\alpha_i\rangle x_\alpha$;

4) $[x_\alpha,x_{-\alpha}]=h_\alpha=$ an integer linear combination of~$h_i$;

5) $[x_\alpha,x_\beta]=N_{\alpha \beta} x_{\alpha+\beta}$, 
if $\alpha+\beta\in \Phi$ ($N_{\alpha\beta}\in \mathbb Z$);

6) $[x_\alpha,x_\beta]=0$, if $\alpha+\beta\ne 0$, $\alpha+\beta\notin \Phi$.

A \emph{representation} of a Lie algebra~$\cal L$  in a linear space~$V$
is a linear mapping
$\pi: {\cal L}\to gl(V)$, where
$$
\forall x,y\in {\cal L}\ \pi ([x,y])=\pi(x)\pi(y)-\pi(y)\pi(x).
$$
A representation is called \emph{faithful} if it has the zero kernel.

\subsection{Elementary Chevalley groups}\leavevmode
                                         
More details about elementary Chevalley groups can be found 
in the book~\cite{Steinberg}.

Let $\cal L$ be a semisimple Lie algebra (over~$\mathbb C$)
with a root system~$\Phi$, $\pi: {\cal L}\to gl(V)$ is its
finitely dimensional faithful representation. Then we can choose a basis
in~$V$ so that all operators $\pi(x_\alpha)^n/n!$ for
$n\in \mathbb N$ are integer (nilpotent) 
matrices. An integer matrix can be considered also as a matrix
over an arbitrary  commutative ring with a unit. Let $R$ be such a ring.
Consider matrices  $n\times n$
over~$R$, the matrices 
$\pi(x_\alpha)^n/n!$ for $\alpha\in \Phi$, $n\in \mathbb N$
we map into $M_n(R)$.

Now let us consider automorphisms of the free module $R^n$ of the form
$$
\exp (tx_\alpha)=x_\alpha(t)=1+t\pi(x_\alpha)+t^2 \pi(x_\alpha)^2/2+\dots+
t^n \pi(x_\alpha)^n/n!+\dots
$$
Since the matrices $\pi(x_\alpha)$ are nilpotent, this series
is finite. The subgroup
of the group $Aut(R^n)$, generated by all automorphisms of the form
 $x_\alpha(t)$, $\alpha\in 
\Phi$, $t\in R$, is called an \emph{elementary Chevalley group} 
(notation: $E_\pi(\Phi,R)$). 

In an elementary Chevalley group we can introduce the following important
elements and subgroups:

--- the subgroup $U$ is generated by all $x_\alpha(t)$, $\alpha \in \Phi^+$, 
$t\in R$;

--- the subgroup $V$ is generated by all $x_\alpha(t)$, $\alpha \in \Phi^-$, 
$t\in R$;

--- $w_\alpha(t)=x_\alpha(t) x_{-\alpha}(-t^{-1})x_\alpha(t)$, $\alpha\in \Phi$,
$t\in R^*$;

--- $h_\alpha (t) = w_\alpha(t) w_\alpha(1)^{-1}$;

--- the subgroup $N$ is generated by all $w_\alpha (t)$, $\alpha \in \Phi$, $t\in R^*$;

--- the subgroup $H$ is generated by all $h_\alpha(t)$, $\alpha \in \Phi$, $t\in R^*$;

--- $B=HU$.

It is known that the group $N$ is a normalizer of~$H$ in the Chevalley group,
the quotient groups $N/H$ is isomorphic to the Weil group $W(\Phi)$,
$U$ is normal in~$B$.

Elementary Chevalley groups are defined even not by a representation of
a Lie algebra, but by a \emph{weight lattice} of this representation.

Let us define this notion.

If $V$ is a representation space of a Lie algebra~$\cal L$
(with a Cartan subalgebra~$\cal H$), then a functional
 $\lambda \in {\cal H}^*$ is called
a \emph{weight} of this representation, if there exists
a nonzero vector $v\in V$ 
such that for any $h\in {\cal H}$
$$
\pi(h) v=\lambda (h)v.
$$
All weights of a given representation (up to addition) generate
a lattice (free Abelian group with  $\mathbb Z$-basis, that is a
 $\mathbb C$-basis in~${\cal H}^*$),
that is called the \emph{weight lattice} $\Lambda_\pi$.

Elementary Chevalley group is completely
defined by a root system~$\Phi$,
a coomutative ring~$R$ with unit and a weight lattice $\Lambda_\pi$.

Among all lattices we mark two:
the lattice corresponded to the adjoint representation 
and generated by all roots (the \emph{adjoint
lattice}~$\Lambda_{ad}$) 
and the lattice generated by all weights of all representations
 (the \emph{universal lattice}~$\Lambda_{sc}$).
For every faithful representation~$\pi$ we have the inclusion
$$
\Lambda_{ad}\subseteq \Lambda_\pi \subseteq \Lambda_{sc}.
$$
Therefore we have the \emph{adjoint} and the \emph{universal}
elementary Chevalley groups.

Every elementary Chevalley group satisfies the following
conditions:

(R1) $\forall \alpha\in \Phi$ $\forall t,u\in R$\quad $x_\alpha(t)x_\alpha(u)=
x_\alpha(t+u)$;

(R2) $\forall \alpha,\beta\in \Phi$ $\forall t,u\in R$\quad
 $\alpha+\beta\ne 0\Rightarrow$
$$
[x_\alpha(t),x_\beta(u)]=x_\alpha(t)x_\beta(u)x_\alpha(-t)x_\beta(-u)=
\prod x_{i\alpha+j\beta} (c_{ij}t^iu^j),
$$ 
where $i,j$ are integral positive numbers, 
the product is taken by all roots $i\alpha+j\beta$, with
some fixed order; $c_{ij}$ are integral numbers not depending of 
$t$ and~$u$;

(R3) $\forall \alpha \in \Phi$ $w_\alpha=w_\alpha(1)$;

(R4) $\forall \alpha,\beta \in \Phi$ $\forall t\in R^*$
$w_\alpha h_\beta(t)w_\alpha^{-1}=h_{w_\alpha (\beta)}(t)$;

(R5) $\forall \alpha,\beta\in \Phi$ $\forall t\in R^*$ 
$w_\alpha x_\beta(t)w_\alpha^{-1}=x_{w_\alpha(\beta)} (ct)$, where $c=c(\alpha,\beta)=
\pm 1$;

(R6) $\forall \alpha,\beta\in \Phi$ $\forall t\in R^*$ $\forall u\in R$
$h_\alpha (t)x_\beta(u)h_\alpha(t)^{-1}=x_\beta(t^{\langle \beta,\alpha
\rangle} u)$.

By $X_\alpha$ we denote the group generated by all
$x_\alpha (t)$ for $t\in R$.

\subsection{Chevalley groups.}\leavevmode

More details about Chevalley groups can be found in the book~\cite{Steinberg},
and also in the papers \cite{Chevalley}, \cite{Artem_dis}, \cite{VavPlotk1} 
(see also references in these papers). 

A set $X\subseteq \mathbb C^n$ is called an \emph{affine
variety}, if $X$ is a set of common zeros in $\mathbb C^n$
of a finite system of polynomials from $\mathbb C[x_1,\dots,x_n]$.

Topology of an affine $n$-space, where the system of closed sets coincides 
with the system of affine varieties, is called  \emph{Zarisski topology}.
A variety is called  \emph{irreducible}, if it can not be
represented as a union of two  proper nonempty closed
subsets.

Let $G$ be an affine variety with some group structure.
If both mappings
\begin{align*} 
m:& G\times G\to G,& m(x,y)&=xy,\\
i:& G\to G,& i(x)&=x^{-1}
\end{align*}
can be expressed as polynomials of their coordinates, then $G$ 
is called an \emph{affine algebraic group}. 
A \emph{linear algebraic group}  is an arbitrary
algebraic subgroup in $M_n(\mathbb C)$ (with usual matrix multiplication).

An algebraic group is called  \emph{connected}, if it is
irreducible as a variety.

Every algebraic group $G$ contains the single greatest connected
solvable normal subgroup. It is called a \emph{radical~$R(G)$}. 
A connected algebraic group with trival radical
is called \emph{semisimple}.

Semisimple linear algebraic groups over $\mathbb C$ (or any other
algebraically closed field~$K$) are exactly elementary Chevalley
group $E_\pi(\Phi,K)$ (see~\cite{Steinberg},~\S\,5).

All such groups are defined in $SL_n(K)$ as  common  zeros of
polynomials of matrix entries $a_{ij}$ with integer coefficients
 (for example, in the case of the root system $A_l$ and 
the universal representation we have 
$n=l+1$ and the single polynomial $\det (a_{ij})-1=0$).
It is clear that multiplying and taking the inverse element are also 
defined by polynomials with integer coefficients. 
Therefore these polynomials can be 
considered as polynomials over an arbitrary commutative ring with a unit.
Let some elementary Chevalley group $E$ over~$\mathbb C$ be defined
in $SL_n(\mathbb C)$ by polynomials $p_1(a_{ij}),\dots, p_m(a_{ij})$.
For a commutative ring~$R$ with a unit consider
the group
$$
G(R)=\{ (a_{ij}\in M_n(R)| \widetilde p_1(a_{ij})=0,\dots
,\widetilde p_m(a_{ij})=0\},
$$
where $\widetilde p_1(\dots),\dots \widetilde p_m(\dots)$ are polynomials with
the same coefficients as $p_1(\dots)$, \dots, $p_m(\dots)$,
but considered over the ring~$R$.

This group is called a \emph{Chevalley group} $G_\pi(\Phi,R)$ of type~$\Phi$
over~$R$, and for every algebraically closed field~$K$ 
it coincides with an elementary Chevalley group.

The subgroup of all diagonal (in the standard basis of weight vectors) matrices
of a Chevalley group $G_\pi(\Phi,R)$ is called a  \emph{standard
maximal torus} of $G_\pi(\Phi,R)$, it is denoted by $T_\pi(\Phi,R)$. 
This group is isomorphic to $Hom(\Lambda_\pi, R^*)$.

Let us denote by $h(\chi)$
the elements of $T_\pi (\Phi,R)$, corresponded to the homomorphism
 $\chi\in Hom (\Lambda(\pi),R^*)$.

In particular, $h_\alpha(u)=h(\chi_{\alpha,u})$ ($u\in R^*$, $\alpha \in \Phi$),
where
$$
\chi_{\alpha,u}: \lambda\mapsto u^{\langle \lambda,\alpha\rangle}\quad
(\lambda\in \Lambda_\pi).
$$

\subsection{Chevalley groups and elementary Chevalley groups.}\leavevmode

The interrelations between Chevalley groups and the corresponding
elementary subgroups constitute one of the major problems in the 
theory of Chevalley groups over rings.
Whereas for an elementary Chevalley group there is a very nice system
of generators $x_\alpha
(\xi)$, $\alpha\in \Phi$, $\xi\in R$, and the relations among these generators
are fairly well understood, noting like that is available for the Chevalley
group itself.

If $R$ is algebraically closed field, then always
$$
G_\pi (\Phi,R)=E_\pi (\Phi,R)
$$
for any representation~$\pi$. This equality is not true even for
the case of fields, that are not algebraically closed.

But if the group $G$ is simply-connected, and the ring $R$ is \emph{semilocal}
(i.e. contains finite number of maximal ideals), then we have
$$
G_{sc}(\Phi,R)=E_{sc}(\Phi,R)
$$
(see \cite{M}, \cite{Abe1}, \cite{St3}, \cite{AS}).

Besides that, it is known, that the groups $G_{sc}$ and $E_{sc}$ coincide, if
the ring~$R$ is euclidean, Dedekind of arithmetic type~\cite{BMS},
\cite{M}, is the polynomial ring with coefficients in a field or
 a principal ideal 
ring~\cite{Su2}, \cite{SK}, \cite{Ko1}, \cite{A3}, \cite{A4}, \cite{Vo},
\cite{Cs}, \cite{GMV}.

Let us show the distinction between Chevalley groups and their 
elementary subgroups in the case when the ring $R$ 
is semilocal, and the Chevalley group is not simply-connected.
In this case $G_\pi (\Phi,R)=E_\pi(\Phi,R)T_\pi(\Phi,R)$ 
(see~\cite{Abe1}, \cite{AS}, \cite{M}), and the elements
 $h(\chi)$ are connected with elementary generators by the formula
\begin{equation}\label{e4}
h(\chi)x_\beta (\xi)h(\chi)^{-1}=x_\beta (\chi(\beta)\xi).
\end{equation}
It is very well-known that the elementary group
 $E_2(R)=E_{sc}(A_1,R)$ is not necessarily normal 
in the special linear group $SL_2(R)=G_{sc}(A_1,R)$
(see~\cite{Cn}, \cite{Sw}, \cite{Su1}).

But if $\Phi$ is an irreducible root system of rank $l\ge 2$,
then $E(\Phi,R)$ is always normal in $G(\Phi,R)$. 
In the case of semilocal rings
from the formula~\eqref{e4} we see that
$$
[G(\Phi,R),G(\Phi,R)]\subseteq E(\Phi,R).
$$
If the ring $R$ also contains~$1/2$, then it is easy to show that
$$
[G(\Phi,R),G(\Phi,R)]=E(\Phi,R)
$$
(we shall do it later).

\section{Easy theorem (theorem 2).}\leavevmode

In this section we  suppose that the ring $R$ is an infinite field
of characteristic $\ne 2$.

We are going to prove the following theorem:

{\bf Theorem.} \emph{Let $G=G_\pi (\Phi,K)$ and $G'=G_{\pi'}(\Phi,K')$
\emph{(}or $E_\pi (\Phi,K)$ and $E_{\pi'}(\Phi,K'))$ are two
 \emph{(}elementary\emph{)}
Chevalley groups over elementary equivalent fields $K$ and~$K'$,
where representations $\pi$ and $\pi'$ have the weight lattices.
Then the groups $G$ and $G'$ are elementarily equivalent.}

In the book \cite{KeisChang}, p.\,395 it was proved that elementary
equivalence is preserved under taking direct products,
therefore

\begin{proposition}\label{p1}
If the semisimple Lie algebra  ${\cal L}={\cal L}_1\oplus \dots \oplus
{\cal L}_k$, where the algebras ${\cal L}_1,\dots,{\cal L}_k$ are simple, 
$R,R'$ are fields  \emph{(}rings\emph{)}, then mutually equivalence of 
\emph{(}elementary\emph{)} 
Chevalley groups
$G_{\pi|_{{\cal L}_1}}({\cal L}_1,R)$ and $G_{\pi|_{{\cal L}_1}}({\cal L}_1,R')$,
\dots,
$G_{\pi|_{{\cal L}_k}}({\cal L}_k,R)$ and $G_{\pi|_{{\cal L}_k}}({\cal L}_k,R')$
implies elementary equivalence of (elementary) Chevalley groups 
$G_{\pi}({\cal L},R)$ and $G_{\pi}({\cal L},R')$.
\end{proposition}

Therefore, we need to prove our theorem for simple Lie albegras.

The following theorem holds for arbitrary commutative rings $R$ and $R'$
with units.

\begin{theorem}\label{easytheor}
If two Chevalley groups $G=G_\pi(\Phi,R)$ and $G'=G_\pi(\Phi,R')$ are constructed
by the same complex Lie algebra of the type~$\Phi$ and by the same
representation~$\pi$, and by elementarily equivalent rings $R$ and~$R'$,
respectively, then $G\equiv G'$.
\end{theorem}

\begin{proof}
As we know from the definition of Chevalley groups,
\begin{align*}
G&= \{ (a_{ij}\in M_n(R)\vert p_1(a_{ij})=p_2(a_{ij})=\dots =p_m(a_{ij})=0\},\\
G'&= \{ (a_{ij}\in M_n(R')\vert p_1(a_{ij})=p_2(a_{ij})=\dots =p_m(a_{ij})=0\},
\end{align*}
where $p_1,p_2,\dots,p_m$ are some well-known polynomials
with integral koefficients,
$n$ is some known natural number.

Suppose that we have some sentence~$\varphi$ of the group language,
 considered in the groups $G$ and~$G'$. Let us translate it
to the sentence $\widetilde \varphi$ of the ring language as
follows:

--- the subformula $\forall g\ \psi (g)$ is translated to the subformula
$$
\forall a_{11}^g,\dots,a_{nn}^g (p_1(a_{11}^g,\dots,a_{nn}^g)=0\land\dots\land
p_m(a_{11}^g,\dots,a_{nn}^g)=0\Rightarrow \widetilde \psi(a_{11}^g,\dots,
a_{nn}^g));
$$
--- the subformula $\exists g\ \psi (g)$ is translated to the subformula
$$
\exists a_{11}^g,\dots,a_{nn}^g (p_1(a_{11}^g,\dots,a_{nn}^g)=0\land\dots\land
p_m(a_{11}^g,\dots,a_{nn}^g)=0\land \widetilde \psi(a_{11}^g,\dots,
a_{nn}^g));
$$
--- the subformula $g=h$ is translated to the subformula
$$
a_{11}^g=a_{11}^h\land \dots \land a_{nn}^g=a_{nn}^h;
$$
--- the subformula $g=h\cdot f$ is translated to the subformula
$$
\bigwedge_{i,j=1}^n \left( a_{ij}^g=\sum_{k=1}^n a_{ik}^h\cdot a_{kj}^f\right).
$$

It is clear that $G(R)\vDash \varphi$ if and only if $R\vDash \widetilde 
\varphi$.

Therefore, if the rings $R$ and $R'$ are elementarily equivalent, then 
for every sentence~$\varphi$  of the group language 
$$
G\vDash \varphi \Leftrightarrow R\vDash \widetilde \varphi \Leftrightarrow
R'\vDash \widetilde \varphi \Leftrightarrow G'\vDash \varphi.
$$
Consequently, $G\equiv G'$.
\end{proof}

The following theorem holds for fields, local and semilocal
rings  $R$ and $R'$ with $1/2$.

\begin{theorem}
If two elementary Chevalley groups $E=E_\pi(R,\Phi)$ and $E'=E_\pi(R',\Phi)$ 
are constructed by the same complex Lie algebra of the type~$\Phi$ and
the same representation~$\pi$, and by elementarily equivalent
semilocal rings $R$ and~$R'$ with~$1/2$, then $E\equiv E'$.
\end{theorem}

It is clear that this theorem follows from the previous one
and the following proposition:

\begin{proposition}
If the ring $R$ is semilocal and has~$1/2$, then the elementary subgroup
$E=E_\pi(R,\Phi)$ is definable in the Chevalley group $G=G_\pi(R,\Phi)$ without
parameters, i.e., there exists a formula $\varphi_{\pi,\Phi}(x)$ of the group
language with one free variable~$x$, that is true in the group~$G$ on the element
$g\in G$ if and only if $g\in E$.
\end{proposition}

The proof of this proposition is contained in the next section, it is 
the corollary of the stronger statement.

\section{Elementary adjoint groups.}\leavevmode

Now we want to prove that if two (elementary) Chevalley groups
are elementarily equivalent, then their root systems and weight 
lattices coincide, the fields are elementarily equivalent.

We suppose exactly two properties of the considered fields:

1) they have characterics $\ne 2$, and in the case of type $G_2$ 
their characterics $\ne 3$;

2) they are infinite.

We need the first supposition, because the main part of the proof
is based on the involutions, and also for the fields of characterisctics~$2$
and~$3$
there are some difficulties with some types of Chevalley groups.

The second supposition does not restrict the generality of the result,
because for the finite fields $K$ and $K'$ (elementary) Chevalley groups
$G(K)$ and $G(K')$ ($E(K)$ and $E(K')$) are finite, i.\,e. 
$$
G(K)\equiv G(K')\Rightarrow G(K)\cong G(K').
$$
Therefore we can refer to the results, proved earlier (for example,
see the book~\cite{Steinberg}),
that show that in this case $\Phi\cong \Phi'$, 
$K\cong K'$, i.\,e. $\Phi\cong \Phi'$, $K\equiv K'$.

In the beginning we show that

\begin{proposition}\label{elemsubgr}
If two Chevalley groups $G$ and $G'$ are elementarily equivalent,
then their elementary subgroups  $E$ and $E'$ are also elementarily 
equivalent.
\end{proposition}

\begin{lemma}
Let $G=G_\pi(R,\Phi)$ be a Chevalley group, $E=E_\pi(R,\Phi)$ be
its elementary subgroup, $R$ be a semilocal ring with~$1/2$ \emph{(}and
with $1/3$, if $\Phi\cong G_2$\emph{)}. Then $E=[G,G]$ and their
exists such a number~$N$,
depending only on~$\Phi$, but not on~$R$ \emph{(}and even not on
the representation~$\pi$\emph{)},
that every element of the group~$E$ is a product of not more than 
$N$~commutators of the group~$G$.
\end{lemma}
\begin{proof}
Let the root system~$\Phi$ have the rank~$l$.

If  $R$ is a field, then for every element $g\in E$ there is 
the \emph{Bruhat decomposition} (see~\cite{Steinberg}) $g=uhwu'$, $u,v\in U$, $h\in H$, 
$w\in W$. If $R$ is a (semi)local ring, then for every element
$g\in E$ tere is  the  \emph{Gauss decomposition} (see~\cite{Abe1}) $g=uhvu'$,
$u,u'\in U$, $v\in V$, $h\in H$. It is known (see~\cite{Steinberg}), 
that elements $u,u',v$ 
can be represented as a products of not more than~$n$ 
(the number of positive roots of the system~$\Phi$) elements~$x_\alpha(t)$, 
and the element~$h$ is a product of not more
than~$l$ elements of the form~$h_\alpha(t)$, that are products
of not more than six elements~$x_\alpha(t)$.

Therefore,
every element of the group $E_\pi({\cal L},K)$
is a product of not more than $6l+3n$ elements $x_\alpha(t)$, where $n$ is
the number of positive roots, depending on~$l$ as follows:

\begin{center}
{\small

\begin{tabular}{|c|c|c|}
\hline
root system & rank & $n$ \\
\hline 
$A_l$ & $l$ & $(l^2+l)/{2} $ \\
\hline
$B_l$ & $l$ & $l^2 $ \\
\hline
$C_l$ & $l$ & $l^2 $ \\
\hline
$D_l$ & $l$ & $l^2-l $ \\
\hline
$E_6$ & $6$ & $36 $ \\
\hline
$E_7$ & $7$ & $63 $ \\
\hline
$E_8$ & $8$ & $120 $ \\
\hline
$F_4$ & $4$ & $24 $ \\
\hline
$G_2$ & $2$ & $6 $ \\
\hline
\end{tabular}
}
\end{center}

\bigskip

Now we need to show that every $x_\alpha(t)$ is a product of
some upper bounded number of commutators.

If in the ring~$R$ there exists an invertible element~$s$ such that
$s^2-1$ is invertible, then $[h_\alpha(s), x_\alpha(t)]=
x_\alpha ((s^2-1)t)$, therefore every $x_\alpha(t)$ 
is a commutator.

If we have systems $A_l$ ($l\ge 2$), $D_l$ ($l\ge 4$), $E_6$, $E_7$
or $E_8$, then, since all roots are conjugate under the action of the Weil group,
we only need to consider
$x_\alpha(t)$ for some root of~$\Phi$. Let it be
the root $\alpha_1+\alpha_2$, where $\alpha_1,\alpha_2$ are
nonorthogonal simple roots.

Then 
$$
[x_{\alpha_1}(t),x_{\alpha_2}(s)]=x_\alpha(\pm ts),
$$
therefore, every $x_\alpha(t)$ is a commutator.

If we have systems $B_l$ ($l\ge 2$), $C_l$ ($l\ge 3$), $F_4$,
then, since all roots of the same length are conjugate under the action 
of the Weil group, we only need to consider arbitrary two roots
of different lengths. In each of these systems there is
the subsystem $B_2$,
therefore, we can consider only this system.

The roots have the form $\pm e_1,\pm e_2$ (short) and
$\pm e_1 \pm e_2$ (long).

Note that $e_1+ e_2= e_1+ e_2$ and no linear combination
of $ e_1$ and $ e_2$ with natural coefficients, different from
$ e_1+ e_2$, can not be a root, therefore
$$
[x_{e_1}(t),x_{e_2}(s)]=x_{e_1+ e_2} (\pm 2ts),
$$
i.\,e., $x_{e_1+ e_2}(t)$ is a commutator (since $1/2\in R$).

Besides, $ e_1=(e_1 -e_2)+e_2$, consequently
$$
[x_{e_1-e_2}(t),x_{e_2}(s)]=x_{e_1}(\pm 2ts) x_{e_1+e_2}(cts^2),
$$
and we see that $x_{e_1}(t)$ is a product of two commutators.

Hence, every $x_\alpha(t)$ is a product of not more than 2 commutators.

Therefore,  every element of an elementary Chevalley group
$E_\pi (\Phi,R)$ is a product of not more than
 $N$ commutators of the group $G_\pi(\Phi,R)$, where the number $N$ 
depends only on the root system~$\Phi$.
\end{proof}

\emph{Proof of Proposition}~2.

Consider the set of sentences
\begin{multline*}
Define_N:= \forall x_1,\dots,x_N, y_1,\dots ,y_N,z_1,\dots,z_N,
t_1,\dots,t_N \\\exists v_1,\dots,v_N,u_1,\dots,u_N\\
(([x_1,y_1]\cdot\dots\cdot [x_N,y_N])\cdot([z_1,t_1]\cdot \dots \cdot
[z_N,t_N])=[u_1,v_1]\cdot \dots \cdot [u_n,v_n]).
\end{multline*}
This sentence states that every element of the commutant of
the group is a product of not more than
 $N$~commutators. We know that for the groups $G$ and $G'$ there exists
(the same, as they are elementarily equivalent) such~$N$, that
the sentence $Define_N$ holds in both groups.
In this case the formula
$$
Commut_N(x):= \exists u_1,\dots , u_N,v_1,\dots v_n (x=[u_1,v_1]\cdot\dots 
\cdot [u_N,v_N])
$$
defines in both groups $G$ and~$G'$ the subgroups $E$ and~$E'$, respectively.
Therefore these groups are elementary equivalent.
$\square$

\medskip

Naturally, if have two elementarily equivalent  elementary Chevalley
groups $E$ and $E'$,
we also have two elementarily equivalent elementary adjoint 
Chevalley groups $E_{ad}$ and $E_{ad}'$,
which are central quotients of th initial groups.

\section{Identification in the classical cases.}\leavevmode

Now we want to identify classical elementary adjoint Chevalley 
groups over fields with some subgroups of $GL_n(K)$.

$\mathbf{A_l}$. It is clear that $G_{sc}(A_l,K)\cong SL_{l+1}(K)$, therefore
\begin{multline*}
E_{sc}(A_l,K)=[G_{sc}(A_l,K),G_{sc}(A_l,K)]\cong [SL_{l+1}(K),SL_{l+1}(K)]=\\
=SL_{l+1}(K),
\end{multline*}
consequently,
$$
E_{ad}(A_l,K)\cong PSL_{l+1}(K).
$$

$\mathbf{C_l}$. Similarly, $G_{sc}(C_l,K)\cong Sp_{2l}(K)$, $E_{sc}(C_l,K)\cong
[Sp_{2l}(K),Sp_{2l}(K)]=Sp_{2l}(K)$, therefore,
$$
E_{ad}(C_l,K)\cong PSp_{2l}(K).
$$

$\mathbf{D_l}$.
For the group $G(D_l,K)$, there exists an intermediate representation~$\pi$, 
for which $G_\pi(D_l,K)\cong SO_{2l}(K)$.
The group $E_\pi (D_l,K)$ is generated by matrices $E+tE_{i,j}-tE_{l+j,l+i}$
and $E+tE_{l+i,l+j}-tE_{l+j,l+i}$, where $1\le i,j\le l$.

From the other side we know that the groups $G_\pi(D_l,K)$
and $E_\pi(D_l,K)$ ``differ by torus'', i.\,e., in their 
Bruhat decompositions  $T$ is changed to~$H$, 
therefore we need to study possible differences between them.
The group $T$ is just the group of all diagonal matrices
 $D=diag[d_1,\dots,d_{2n}]$,
with the condition $DQD^T=Q$, i.\,e. $T$ consists of matrices 
$$
diag[d_1,\dots,d_l,1/d_l,\dots, 1/d_1].
$$
The group $H$ can be found with the help of elements
$h_\alpha(t)$.
These elements have the form:
\begin{align*}
&diag[t_1,1/t_1,1,\dots, 1, t_1, 1/t_1],\\
& diag [s_1,s_1,1,\dots,1,1/s_1,1/s_1],\\
&diag [1,t_2,1/t_2,1,\dots,1,t_2,1/t_2],\\
&diag [1,s_2,s_2,1,\dots,1, 1/s_2,1/s_2,1],\\
&\dots, \\
& diag[1,\dots,1,t_{n-1},1/t_{n-1},t_{n-1}, 1/t_{n-1},1,\dots,1],\\
&diag[1,\dots,1, s_{n-1},s_{n-1},1/s_{n-1},1/s_{n-1},1,\dots,1].
\end{align*}
Note that we can generate $diag[d_1,1,\dots,1,1/d_1]$ only as
$$
diag
[t_1s_1,s_1/t_1,1,\dots,1,t_1/s_1,1/(t_1s_1)],
$$
therefore $d_1=t_1^2$, 
i.\,e., everything depends of squares in the field~$K$.

We get that the group $H$ consists of all diagonal matrices 
$$
diag[d_1,\dots,d_l, 1/d_l,\dots,1/d_1],
$$
where  $(d_1\dots d_l)$ is a square in~$K$, 
consequently, involutions of $E_\pi(D_l,K)$ coincide with involutions
of $G_\pi(D_l,K)$ for $i\in K$, and for $i\notin K$ we have not involutions
with odd number of~$-1$
among the first~$l$ elements of the diagonal.

To get our group $E_{ad}(D_l,K)$, we need to factor the group
 $E_\pi(D_l,K)$ by its center that is trivial for $i\notin K$ and
odd~$l$, and consists of $\pm E_{2l}$ in other cases.

$\mathbf{B_l}$.
Similarly to the case $D_l$, in this case the group $E_{ad}(B_l,K)$
coincides with $SO_{2l+1}(K)$, if every element of~$K$ is
a square, and differs by such a part of torus, that consists
of elements
$$
diag[d_1,\dots,d_l,1,1/d_l,\dots,1/d_1],
$$
where $(d_1\dots d_l)$ is not a square.

Respectively, for $i\notin K$ this groups does not contain
involutions with odd number of $-1$ among the first $l$ elements.

\section{Study of involutions for classical Chevalley groups.}\leavevmode

\subsection{Involutions in the group $PSL_n(K)$.}\leavevmode

Naturally, every involution in the group $PSL_n(K)$ is either the image
of the involution from $SL_n(K)$ (\emph{involution of the first type}), 
or is the image
of such a matrix $A\in SL_n(K)$, that $A^2=\lambda E$ (\emph{involution
of the second type}).

At first  consider involutions of the first type.
Involutions from $SL_n(K)$ 
have in some basis the form $diag[\underbrace{-1,\dots,-1}_{k},\underbrace{1,
\dots,1}_{n-k}]$, $k$ is even. For even $n$ the involutions
 $diag[\underbrace{-1,\dots,-1}_{k},
\underbrace{1,
\dots,1}_{n-k}]$ and $diag[\underbrace{-1,\dots,-1}_{n-k},\underbrace{1,
\dots,1}_{k}]$ are equal in the group $PSL_n(K)$.

For every even~$k$ any two involutions with the same~$k$ are conjugate,
i.e. for even~$n$
there exist $[n/4]$ conjugacy classes of first type involutions,
and for odd~$n$ 
there are $[n/2]$ conjugacy classes.

The centralizer of a first type involution (for $2k\ne n$) consists
of blocks
$$
\begin{pmatrix}
A& 0\\
0& B
\end{pmatrix},\quad A\in GL_k(K),\ B\in GL_{n-k}(K),\ det A\cdot det B=1.
$$
Its center consists of the matrices
$$
diag[\underbrace{a,\dots,a}_{k}\underbrace{b,\dots, b}_{n-k}],\quad a^k\cdot b^{n-k}
=1.
$$

Let us show that there are infinitely many such matrices.

If $\text{G.C.D.}(k,n-k)\ne 1$, then let us change $k$ and $n-k$ to
 $k'=k/\text{G.C.D.}(k,n-k)$
and $l=(n-k)'=(n-k)/\text{G.C.D.}(k,n-k)$. Then we need to prove
that there are infinitely many solutions of the equation $a^{k'}b^l=1$.
For every $\xi\in K^*$ 
elements $a=\xi^{ml}$ and $b=\frac{1}{\xi^{mk'}}$ for $m>0$ satisfy
the condition, therefore
we need to prove that there are infinitely many pairs $(a,b)$.
For $char \, K=0$ it is clear.
Let $char \,K\ne 0$. In any case either $l$, or $k$ is mutually simple with
$p=char\, K$. Let it be $l$. It is clear that, varying $\xi$ and $m$ 
(we know that there are infinitely many such $\xi$), 
we get infinitely many different~$a$.

Therefore, if $2k\ne n$, then the center of the centralizer of
the first type involution is necessary infinite. Its commutant consists
of matrices
$$
\begin{pmatrix}
A& 0\\
0&B
\end{pmatrix},\quad A\in SL_k(K), B\in SL_{n-k}(K),
$$
factorized by the center, and the central 
quotient of this commutant is
$$
PSL_k(K)\times PSL_{n-k}(K).
$$

Let now $2k=n$, i.e., our involution in some basis has the form $diag[
-E_k,E_k]$. Naturally, in this case $k$ is even.

The centralizer of such an involution consists of matrices 
$$
\begin{pmatrix}
A& 0\\
0& B
\end{pmatrix} \quad A,B\in GL_k(K), det\, A\cdot det\, B=1,
$$
and also
$$
\begin{pmatrix} 
0& C\\
D& 0
\end{pmatrix},\quad C,D\in GL_k(K),\ det\, C\cdot det\, D=1.
$$
The center of this centralizer is finite and has the order 2.
It commutant consists of matrices
$$
\begin{pmatrix}
A& 0\\
0& B
\end{pmatrix} \quad A,B\in SL_k(K),
$$
and its central quotient is 
$$
PSL_k(K)\times PSL_k(K).
$$

Now let us go to the second type involutions.

Let us have some semi-involution in the group $SL_n(K)$, i.e. 
a matrix~$A$, such that
$A^2=\lambda E$. Since $char\, K\ne 2$,  the matrix~$A$
is diagonalizable (in~$\overline K$), therefore its proper
values are equal to $\pm \sqrt \lambda$.
So we have $\lambda^n=1$.

Here there are different cases:

1) $\sqrt \lambda\in K$ and $\sqrt \lambda^n=1$. Such a semi-involution 
coincides with a usual involution in $PSL_n(K)$.

2) $\sqrt \lambda\in K$ and $\sqrt \lambda^n=-1$. If $n$  
is odd, then we can put
$\lambda'=\sqrt \lambda$ and go to the case~1). 
Therefore we can suppose that $n$ 
is even. In this case we have involutions of the form
$$
diag [\underbrace{-\sqrt\lambda,\dots,-\sqrt \lambda}_{k},
\underbrace{\sqrt \lambda,\dots,
\sqrt \lambda}_{n-k}],
$$
$k$ is odd. There are more $[n/4]$ conjugate classes of these involutions.

The centralizer of such a $(k,n-k)$-involution consists of matrices
$$
\begin{pmatrix}
A& 0\\
0& B
\end{pmatrix} \quad A\in GL_k(K), B\in GL_{n-k}(K), det\, A\cdot det\, B=1,
$$
and for the odd number $k=n/2$ we also have matrices
$$
\begin{pmatrix}
0& C\\
D& 0
\end{pmatrix} \quad C,D\in GL_k(K), det\, C\cdot det\, D=-1.
$$
In all cases, except the last one, the center of
the centralizer is infinite, the central quotient of the commutant 
of the centralizer is
$$
PSL_{k}(K)\times PSL_{n-k}(K),
$$
if $k=1$, then it is just $PSL_{n-1}(K)$.

3) $\sqrt \lambda\notin K$, but $\lambda\in K$. Note that
there are equal number of the proper values 
$+\sqrt \lambda$ and $-\sqrt \lambda$, since
$tr\, A\in K$. We have $(-\lambda)^{n/2}=1$.

The matrix 
$$
\begin{pmatrix}
0& \lambda& &        &      &  &\\
1& 0      & &        &      &  &\\
 &        &0& \lambda&      &  &\\
 &        &1&     0  &      &  &\\
 &        & &        &\ddots&  &\\
 &        & &        &      &0&\lambda\\
 &       &  &        &      &1& 0
\end{pmatrix}
$$
represents such an involution.
We also can write it in the form
$$
\begin{pmatrix}
0& \lambda E\\
E& 0
\end{pmatrix}.
$$

\subsection{Study of involutions in groups of the type $C_l$.}\leavevmode

Since the center of $Sp_{2l}(K)$ consists of two elements $\pm E$, we have that
involutions of
$PSp_{2l}(K)$ are such elements~$\widetilde A$ 
that $A^2=E$ (the first type) and 
such elements~$\widetilde A$ that $A^2=-E$ (the second type).

At first consider the first type involutions: in some basis such an
involution $A$ is
$$
\begin{pmatrix}
E& 0\\
0& -E
\end{pmatrix},
$$
in this basis the form $Q$ is
$$
Q_A=\begin{pmatrix}
Q_1& Q_2\\
-Q_2^T& Q_3
\end{pmatrix}
$$
Then the condition $AQ_AA^T=Q_A$ implies
\begin{multline*}
\begin{pmatrix}
E& 0\\
0& -E\end{pmatrix}
\begin{pmatrix}
Q_1& Q_2\\
-Q_2^T& Q_3
\end{pmatrix}\begin{pmatrix}
E& 0\\
0& -E
\end{pmatrix}=\begin{pmatrix}
Q_1& Q_2\\
-Q_2^T& Q_3
\end{pmatrix}\Rightarrow\\
\begin{pmatrix}
Q_1& -Q_2\\
Q_2^T& Q_3
\end{pmatrix}
=\begin{pmatrix}
Q_1& Q_2\\
-Q_2^T& Q_3
\end{pmatrix}\Rightarrow
Q_2=0.
\end{multline*}
Therefore $Q_1$ and $Q_3$ are non-degenerate skew-symmetric matrices, i.e. 
have even dimension.

It is clear, that for every even $k\le l$ $(2k,2(l-k))$-involutions of
the first type exist in the group $PSp_{2l}(K)$, since we can take
the matrices
$$
diag[ \underbrace{-1,\dots,-1}_k, \underbrace{1,\dots,1}_{2(l-k)},
\underbrace{-1,\dots,-1}_k].
$$

The centralizer of such an involution consists of matrices
$$
\begin{pmatrix}
A& 0\\
0& B
\end{pmatrix}, \quad A\in Sp_{2k}(K), B\in SP_{2(l-k)}(K),\text{ or } 
A\in Sp_{2k}^-(K), B\in Sp_{2(l-k)}^-(K)
$$
for $2k\ne l$.

Its center is finite and consists of two elements.

The commutant of the centralizer consists of the matrices
$$
\begin{pmatrix}
A& 0\\
0& B
\end{pmatrix} \quad A\in Sp_{2k}(K),B\in Sp_{2(l-k)}(K),
$$
   and the central quotient of this commutant is
$$
PSp_{2k}(K)\times PSp_{2(l-k)}(K).
$$
The exception can be in the case of small~$k$: for 
$k=1$ 
$$
Sp_2(K)\cong SL_2(K),
$$
i.\,e., we have the group
$$
PSL_2(K)\times PSp_{2l-2}(K).
$$

Now let us look what happens for $2k=l$.

Let $Q_{2l}=diag[Q_l,Q_l]$, $I=diag[E_l,-E_l]$.
The condition
$$
\begin{pmatrix}
A& B\\
C& D
\end{pmatrix}
\begin{pmatrix}
E_l& 0\\
0& -E_l
\end{pmatrix}=
\begin{pmatrix}
-E_l& 0\\
0& E_l
\end{pmatrix}
\begin{pmatrix}
A& B\\
C& D
\end{pmatrix}
$$
can appear, therefore $A=D=0$, i.e. the centralizer consists of
matrices
$$
\begin{pmatrix}
A& 0\\
0& D
\end{pmatrix},\quad A,D\in Sp_{l=2k}(K)
$$
and
$$
\begin{pmatrix}
0& B\\
C& 0
\end{pmatrix},\quad B,C,\in Sp_{2k}(K).
$$
Its center consists of two elements, and the commutant 
consists of matrices
$$
\begin{pmatrix}
A& 0\\
0& D
\end{pmatrix},\quad A,D\in Sp_{2k}(K).
$$

Now consider the second type involutions. Since $A^2=-E$, then in the field
 $\overline K$ the involution
$A$ is diagonalizable with proper values $\pm i$:
$$
\widetilde A=\begin{pmatrix}
iE_k& 0\\
0& -iE_{2l-k}
\end{pmatrix},
$$
and the form $Q$ has in this basis the form
$$
Q_A=\begin{pmatrix}
Q_1& Q_2\\
-Q_2^T& Q_3
\end{pmatrix}.
$$
The condition
$$
\begin{pmatrix}
iE_k& 0\\
0& -iE_{2l-k}
\end{pmatrix}
\begin{pmatrix}
Q_1& Q_2\\
-Q_2^T& Q_3
\end{pmatrix}
\begin{pmatrix}
iE_k& 0\\
0& -i E_{2l-k}
\end{pmatrix}=
\begin{pmatrix}
Q_1& Q_2\\
-Q_2^T& Q_3
\end{pmatrix}
$$
implies $Q_1=Q_3=0$, therefore $k=l$.

Thus,
$$
\widetilde A=\begin{pmatrix}
iE_l& 0\\
0& -iE_l
\end{pmatrix},\
Q_A=\begin{pmatrix} 0& Q_2\\
-Q_2^T& 0
\end{pmatrix}.
$$
In both cases (if  $i\in K$, or not) such involutions are contained 
in $PSp_{2l}(K)$:
for example,
$$
A=Q=\begin{pmatrix}
0& E\\
-E& 0
\end{pmatrix}.
$$

To find the centralizer of~$A$, let us consider 
(if needed, then in~$\overline K$)
the conjugate involution
$$
I=\begin{pmatrix}
iE& 0\\
0& -iE
\end{pmatrix},
$$
the form $Q$ is
$$
\begin{pmatrix}
0& E\\
-E& 0
\end{pmatrix}.
$$

Let us have the matrix
$$
M=\begin{pmatrix}
A& B\\
C& D
\end{pmatrix}.
$$
Then there are two possibilities:

1) $MI=IM$ in $Sp_{2l}(K)$,
therefore $B=C=0$, consequently, $AD^T=E$.
Thus,
$$
M=\begin{pmatrix}
A&0\\
0& (A^T)^{-1}
\end{pmatrix},\quad A\in GL_l(K).
$$

2) $MI=-IM$ in $Sp_{2l}(K)$,
therefore $A=D=0$, consequently, $BC^T=-E$.
Thus,
$$
M=\begin{pmatrix}
0&B\\
-(B^T)^{-1}& 0
\end{pmatrix},\quad B\in GL_l(K).
$$

If $l$ is odd, then the matrices of the second form
have determinant~$1$, i.e., they are contained  in our 
group, therefore, the center of the centralizer in this case
consists of two elements.
If $l$ is even, then the matrices of the second form have
determinant~$-1$, i.e.,
they are not contained in our group, therefore the center of the centralizer is
infinite.

The commutant in the both cases has the form
$$
\begin{pmatrix}
A& 0\\
0& (A^T)^{-1}\end{pmatrix},\quad A\in SL_l(K).
$$
The central quotient of this commutant is isomorphic to $PSL_l(K)$.

\subsection{Study of involutions in groups of
the type $B_l$.}\leavevmode

As we have already shown, for the case $E_{ad}(B_l,K)$ 
the situations $i\in K$ and $i\notin K$ are fundamentally different. 
Let us consider them separately.

1) $i\in K$. Since $E_{ad}(B_l,K)$ is a subgroup in $SO_{2l+1}(K)$, 
we have that all involutions in $E_{ad}(B_l,K)$ are the first type
involutions.
If $i\in K$, then all involutions of $SO_{2l+1}(K)$ are contained
in $E_{ad}(B_l,K)$, therefore we need to study namely these 
involutions.
It is clear, that involutions of $SO_{2l+1}(K)$
have in some basis the form
$diag[\underbrace{-1,\dots,-1}_{2k},\underbrace{1,\dots,1}_{2l-2k+1}]$, 
there are exactly  $l$~conjugate classes of them.
The centralizer of such an involution consists of matrices
$$
\begin{pmatrix}
A& 0\\
0& B
\end{pmatrix},\quad A\in EO_{2k}(K), B\in EO_{2l-2k+1}(K)
$$
or
 $$
A\in O^-_{2k}(K), B\in O_{2l-2k+1}^-(K).
$$

If $2k=2l$, then we have the matrix
$$
\begin{pmatrix}
A& 0\\
0& \pm 1
\end{pmatrix},\quad A\in O_{2l}(K).
$$

If $k=1$, then we have the matrices
$$
\begin{pmatrix}
a& 0&0 \\
0& 1/a&0 \\
0&0 & B
\end{pmatrix},\quad B\in SO_{2l-1}(K)
$$
and
$$
\begin{pmatrix}
0& a&0 \\
1/a& 0&0 \\
0&0 & B
\end{pmatrix},\quad B\in O_{2l-1}^-(K).
$$

The center of this centralizer has the order~$2$, like all
other centers of centralizers. The central quotients
of commutants of these centralizers are isomorphic
to
$$
EO_{2l-1}(K), PEO_{4}(K)\times EO_{2l-3}(K).\dots,
PEO_{2l-2}(K)\times EO_3(K), PEO_{2l}(K).
$$

2) $i\notin K$. As we have shown,  
the involutions of $SO_{2l+1}(K)$, 
that have $2(2i+1)$ elements $-1$ on the diagonal, are not contained in
$E_{ad}(B_l,K)$. Therefore, in this group we have the involutions
$$
diag[-1,-1,-1,-1,1,\dots,1], \dots, 
diag[\underbrace{-1,\dots,-1}_{4k},1,\dots,1],
\dots.
$$
It is clear that centers of their centralizers also consist of
two elements, and 
the central quotients of commutants of these centralizers are
$$
PEO_4(K)\times EO_{2l-3}(K),\dots, PEO_{2l-2}(K)\times EO_3(K),\text{ or }
PEO_{2l}(K).
$$

\subsection{Study of involutions for groups of the type $D_l$ $(l\ge 4)$.}\leavevmode

In this case, like in the previous one, 
the situations $i\in K$ and $i\notin K$ are different.

1) $i\in K$. At first we shall consider the first type involutions.
These are  involutions of $PSO_{2l}(K)$  with
 inverse images being involutions of $SO_{2l}(K)$. If $i\in K$,
then all involutions of $SO_{2l}(K)$ are contained in $EO_{2l}(K)$,
therefore all first type involutions of $PSO_{2l}(K)$
are contained in $PEO_{2l}(K)$.
It is clear, that involutions of $SO_{2l}(K)$ in some basis have the form
$diag[\underbrace{-1,\dots,-1}_{2k}, \underbrace{1,\dots,1}_{2(l-k)}]$,
there are  $[l/2]$ conjugate classes of them.

For $2k=l$ the centralizer of such an involution consists of matrices 
$$
\begin{pmatrix}
A& 0\\
0& B
\end{pmatrix},\quad A\in EO_{2k}(K),\ B\in EO_{2(l-k)}(K).
$$
Its center consists of 2 elements: $E$ and the involution itself.
The central quotient of the commutant of this centralizer is
$$
PEO_{2k}(K)\times PEO_{2(l-k)}(K).
$$
The exception can be for small~$k$: for $k=1$
$$
PEO_2(K)\times PEO_{2(l-1)}(K)\cong PEO_{2(l-1)}(K).
$$
Note that the type $D_2$ coincides with $A_1\times A_1$, and the type $D_3$ 
coincides with $A_3$.

Now let us look for $2k=l$. Is is clear that $k$~is even.
Let 
$$
Q_{2l}=\begin{pmatrix}
Q_l& 0\\
0& Q_l
\end{pmatrix},\quad I=\begin{pmatrix}
E_l& 0\\
0& -E_l
\end{pmatrix}.
$$

Then it is possible
$$
\begin{pmatrix}
A& B\\
C& D
\end{pmatrix}
\begin{pmatrix}
E_l& 0\\
0& -E_l
\end{pmatrix}=
\begin{pmatrix}
-E_l& 0\\
0& E_l
\end{pmatrix}
\begin{pmatrix}
A& B\\
C& D
\end{pmatrix}\Rightarrow A=D=0,
$$
i.e. the whole cntralizer consists of matrices
\begin{gather*}
\begin{pmatrix}
A& 0\\
0& D
\end{pmatrix},\quad A,D\in SO_k(K),\text{ or } A,D\in O^-(K)\\
\begin{pmatrix}
0& B\\
C& 0
\end{pmatrix}\quad B,C\in SO_k(K),\text{ or } B,C\in O^-_k(K).
\end{gather*}
Its center consists of two elements, and the commutant consists of
matrices
$$
\begin{pmatrix}
A& 0\\
0& D
\end{pmatrix},\quad A,D\in EO_k(K).
$$
Now let us look for the second type involutions.
Since $A^2=-E$, we have that in the field $K$ ($i\in K$) the involution
$A$ is diagonalizable with proper values
$\pm i$, the matrix $A$ is conjugate to
$$
\widetilde A=\begin{pmatrix}
iE_k& 0\\
0& -iE_{2l-k}
\end{pmatrix},
$$
and the form~$Q$ in this conjugate basis is
$$
Q_A=\begin{pmatrix}
Q_1& Q_2\\
Q_2^T& Q_3
\end{pmatrix}.
$$
We have
$$
\begin{pmatrix}
iE_k& 0\\
0& -i E_{2l-k}
\end{pmatrix}
\begin{pmatrix}
Q_1& Q_2\\
Q_2^T& Q_3
\end{pmatrix}
\begin{pmatrix}
iE_k& 0\\
0& -iE_{2l-k}
\end{pmatrix}=
\begin{pmatrix}
Q_1& Q_2\\
Q_2^T& Q_3
\end{pmatrix},
$$
therefore $Q_1=Q_3=0$, i.e. $k=l$.
Consequently,
$$
\widetilde A=\begin{pmatrix}
iE_l& 0\\
0& -iE_l
\end{pmatrix},\quad Q_A=\begin{pmatrix}
0& Q_2\\
Q_2^T& 0
\end{pmatrix}.
$$
For example,
$$
Q_A=\begin{pmatrix} 0& E\\
E& 0
\end{pmatrix}.
$$

The centralizer of $\widetilde A$ consists of two sets:

1) 
\begin{multline*}
\begin{pmatrix}
A& B\\
C& D
\end{pmatrix}
\begin{pmatrix}
iE& 0\\
0& -iE
\end{pmatrix}
=\begin{pmatrix}
iE& 0\\
0& -iE
\end{pmatrix}
\begin{pmatrix}
A& B\\
C& D
\end{pmatrix}\Rightarrow B=C=0\Rightarrow\\
\Rightarrow  AD^T=E\Rightarrow D=(A^T)^{-1}, A\in GL_l(K).
\end{multline*}

2)
\begin{multline*}
\begin{pmatrix}
A& B\\
C& D
\end{pmatrix}
\begin{pmatrix}
iE& 0\\
0& -iE
\end{pmatrix}
=-\begin{pmatrix}
iE& 0\\
0& -iE
\end{pmatrix}
\begin{pmatrix}
A& B\\
C& D
\end{pmatrix}\Rightarrow A=D=0\Rightarrow\\
\Rightarrow  CB^T=E\Rightarrow C=(B^T)^{-1}, B\in GL_l(K).
\end{multline*}

If $l$ is even, then the matrices of the second form have
determinant~$1$, so they are contained in our group,
i.\,e., the center of this centralizer consists of two elements. If $l$
is odd, then the matrices of the second form have determinant $-1$, i.\,e.,
they are not contained in our group,
therefore in this case the center of our centralizer is infinite.

The central quotient of the commutant of the centralizer in any case is
the group $PSL_l(K)$.

2) $i\notin K$.
As we remember, in this case the number of non-conjugate first type
involutions in the group is smaller. Namely, the group $E_{ad}(D_l,K)$ 
does not contain such involutions of $SO_{2l}(K)$, that have $2(2j+1)$ 
elements $-1$ on the diagonal.
Therefore, in this group we have involutions
$$
diag[-1,-1,-1,-1,1,\dots,1],\dots, diag[\underbrace{-1,\dots,-1}_{4k},1,\dots,1],
\dots
$$
It is clear that the centers of the centralizers also consist of two elements,
and the central quotients of commutants of the centralizers
are
\begin{multline*}
PEO_4(K)\times PEO_{2l-4}(K),\dots ,\\
PEO_{2l-2}(K)\times PEO_2(K),\text{ or } PEO_{2l-4}(K)\times PEO_4(K).
\end{multline*}
In the case of odd $l$ all first type involutions are contained in the group,
and in 
the case of even~$l$ only  half of classes of involutions are contained in it.

If  $l$ is odd, then there is no second type involution
conjugate to
$$
\begin{pmatrix}
iE_l& 0\\
0& -iE_l
\end{pmatrix}.
$$
If $l$ is even, then there exists such an involution, for example, it is
$$
A=\begin{pmatrix}
0& 1& & & \\
-1& 0& & & \\
& &\ddots& &\\
& &      &0& 1\\
& &      &-1& 0
\end{pmatrix}
$$
for 
$$
Q_A=\begin{pmatrix}
0& E_l\\
E_l& 0
\end{pmatrix}.
$$

\section{Formulas separating different classical Chevalley groups.}\leavevmode

In this section we show that for any two classical Chevalley groups
with nonisomorphic root systems there exists a first order sentence
that holds in the first group and does not hold in the second one.

\begin{lemma}\label{klass1}
There exists a first order sentence $\varphi_{A_1}$ that holds in
any adjoint groups of the type $A_1$ and does not hold in all
adjoint Chevalley groups of other classical types.
\end{lemma}

\begin{proof}
In the group $PSL_2(K)$ there is just one conjugate class of involutions:
these are involutions conjugate to
$$
\begin{pmatrix}
0& 1\\
-1& 0
\end{pmatrix}.
$$

The centralizer of such an involution consists of matrices
$$
\begin{pmatrix}
a& b\\
-b& a
\end{pmatrix},\quad a^2+b^2=1,
$$
and
$$
\begin{pmatrix}
a& b\\
b& -a
\end{pmatrix},\quad -a^2-b^2=1.
$$
The commutant of this centralizer consists of
$$
\begin{pmatrix}
a& b\\
-b& a
\end{pmatrix},\quad a^2+b^2=1,
$$
i.\,e., it is commutative.

Let us look if there exist other Chevalley groups, satisfying this condition:
since there is only one conjugate class of involutions, we need
to consider only small dimensions.

$A_2$ is the group $PSL_3(K)$, it has the involution
$$
\begin{pmatrix}
-1& 0&0\\
0& -1& 0\\
0& 0& 1
\end{pmatrix},
$$
its centralizer is isomorphic to $GL_2(K)$, i.\,e., is does not
satisfy the condition that
its commutant is commutative.

$A_3$ is the group $PSL_4(K)$, it has at least the involution $diag[
-1,-1,1,1]$, the commutant of its centralizer is the central product
of $SL_2(K)$ and $SL_2(K)$. It is clear, that it is not commutative.

In the groups $A_l$ for $l\ge 4$ there are at least two conjugate 
classes of involutions, therefore we can not to consider these cases.

$C_3$ is the group $PSp_6(K)$. For the involution
$diag[-1,-1,1,1,1,1]$ the commutant of its centralizer
can not be commutative.

For the groups $C_l$ for $l\ge 4$ there are at least two 
non-conjugate involutions,
therefore we can omit these types.

$B_2$ is the group $EO_5(K)$, it contains at least the involution
$$diag[-1,-1,-1,-1,1],$$ 
with the commutant of centralizer  $EO_4(K)$,
this group is not commutative. 
It is clear that for the groups $B_l$, $l\ge 3$, this commutant can not be 
commutative.

$D_4$ is the group $PEO_8(K)$, it contains the involution
$diag [-1,-1,-1,-1,1,1,1,1]$, the commutant of its centralizer 
 also is not commutative.

Therefore we can write a sentence that separate the group $PSL_2(K)$
from all other classical groups:
\begin{multline*}
\varphi_{A_1}: = \forall M_1\forall M_2 (M_1^2=M^2=1\land M_1\ne 1\land M_2\ne 1
\Rightarrow\\
\Rightarrow \exists X(XM_1X^{-1}=M_2))\land \forall X_1\forall X_2 (\exists Y_1\exists Y_2
\exists Z_1\exists Z_2\exists M (M^2=1\land\\
M\ne 1\land Y_1M=MY_1\land Y_2M=MY_2\land Z_1M=MZ_1\land Z_2M=MZ_2\land\\
X_1=Y_1Z_1Y_1^{-1}Z_1^{-1}\land X_2=Y_2Z_2Y_2^{-1}Z_2^{-1})
\Rightarrow X_1 X_2=X_2X_1).
\end{multline*}
\end{proof}

Note that if for some Chevalley a set $\cal N$ is definable (and it is classical Chevalley group
over~$K$ itself),
then we can write a formula  stating that it is 
isomorphic to the group
$PSL_2(K)$.

It is the following formula:
\begin{multline*}
\varphi_{A_1}^{\cal N}: = \forall M_1\in {\cal N}\forall M_2\in {\cal N} 
(M_1^2=M^2=1\land M_1\ne 1\land M_2\ne 1
\Rightarrow\\
\Rightarrow \exists X\in {\cal N}(XM_1X^{-1}=M_2))\land\\
\land \forall X_1\in {\cal N}
\forall X_2\in {\cal N} (\exists Y_1\in {\cal N}\exists Y_2\in {\cal N}
\exists Z_1\in {\cal N}\exists Z_2\in {\cal N}\exists M\in {\cal N}
 (M^2=1\land\\
M\ne 1\land Y_1M=MY_1\land Y_2M=MY_2\land Z_1M=MZ_1\land Z_2M=MZ_2\land\\
X_1=Y_1Z_1Y_1^{-1}Z_1^{-1}\land X_2=Y_2Z_2Y_2^{-1}Z_2^{-1})\Rightarrow\\
\Rightarrow X_1 X_2=X_2X_1).
\end{multline*}

\begin{lemma}\label{klass2}
There exists a first order sentence $\varphi_{A_2}$ that holds in
any adjoint groups of the type $A_2$ and does not hold in all
adjoint Chevalley groups of other classical types.
\end{lemma}

\begin{proof}
This group is characterized by follows:

1) it has only one conjugate class of  involutions;

2) the centralizer of any involution has an infinite center;

3) the central quotient of the commutant of the centralizer of
any involution is  $PSL_2(K)$, 
i.e. a group satisfying the sentence $\varphi_{A_1}$.

It is clear that these three conditions completely characterize the group
 $PSL_3(K)$.
\end{proof}

\begin{lemma}\label{klass3}
There exists a first order sentence $\varphi_{A_3}$ that holds in
any adjoint groups of the type $A_3$ and does not hold in all
adjoint Chevalley groups of other classical types, except the type~$B_2$.
\end{lemma}
\begin{proof}
We have the group $PSL_4(K)$.

In this case the number of conjugate classes of involution depends
of the basic field, but we know that there exists an involution
with central quotient of the commutant of its centralizer isomorphic
to  $PSL_2(K)\times PSL_2(K)$, 
i.\,e. it is such a group that there exist
matrices $X_1,X_2,Y_1,Y_2$
(for example,
$$
\left.\begin{pmatrix}
1& 1& 0& 0\\
0& 1& 0& 0\\
0& 0& 1& 0\\
0& 0& 0& 1
\end{pmatrix},\
\begin{pmatrix}
2& 0& 0& 0\\
0& 1/2& 0& 0\\
0& 0& 1& 0\\
0& 0& 0& 1
\end{pmatrix},\
\begin{pmatrix}
1& 0& 0& 0\\
0& 1& 0& 0\\
0& 0& 1& 1\\
0& 0& 0& 1
\end{pmatrix},\
\begin{pmatrix}
1& 0& 0& 0\\
0& 1& 0& 0\\
0& 0& 2& 0\\
0& 0& 0& 1/2
\end{pmatrix}\right),
$$
that 
$X_1Y_1=Y_1X_1$, $X_1Y_2=Y_2X_1$, $X_2Y_1=Y_1X_2$,
$X_2Y_2=Y_2X_2$, the centralizer of the set $\{ X_1,X_2\}$
satisfies the formula $\varphi_{A_1}$,
the centralizer of the set $\{ Y_1,Y_2\}$ satisfies the formula
 $\varphi_{A_1}$,
these centralizers commute, any element of the central quotient
of the commutant of the centralizer of our involution is represented
as a product of an element from the first centralizer and an element from the second
one.

It is clear that there are no groups of the type $A_l$  
with the central quotient of the commutant of the centrilizer of some 
involution elementary equivalent to $PSL_2(K)\times PSL_2(K)$.

In the case of the type $C_l$ the group of the type $C_3$ (i.\,e.,
 $PSp_6(K)$)
the involution $diag[-1,-1,1,\dots,1]$ has the central quotient
of the commutant of its centralizer isomorphic to
 $PSL_2(K)\times PSp_4(K)$, it can not be elementary equivalent
to  $PSL_2(K)\times 
PSL_2(K)$, and the involution $diag[i,i,i,-i,-i,-i]$ has the central quotient
of commutant of its centralizer isomorphic to $PSL_3(K)$, it also can not
be elementary equivalent to
$PSL_2(K)\times PSL_2(K)$.

It is clear that the cases $l> 3$ can be omitted, 
since there are products of groups of greater dimensions in these cases.

In the cases $D_l$ we only need to consider groups of the type $D_4$, it is
$PEO_8(K)$. 

The involution $diag [-1,-1,1,\dots,1]$ of this group 
has the central quotient of the commutant of its centralizer
isomorphic to $PEO_6(K)$, it can not be elementary equivalent
to $PSL_2(K)\times PSL_2(K)$.
The involution $diag[-1,-1,-1,-1,1,1,1,1]$ has the central quotient of
the commutant of its centralizer isomorphic to $PEO_4(K)\times PEO_4(K)$, 
it also can not be elementary equivalent to
$PSL_2(K)\times PSL_2(K)$. The involution $diag[i,i,i,i,-i,-i,-i,-i]$ has
the central quotient of the commutant of its centralizer isomorphic to
$PSL_4(K)$, and it is also not our case.

It is clear that for groups of the type $B_l$, $l\ge 3$, i.\,e. 
for groups $EO_{2l+1}(K)$, there is no involution with
central quotient of the commutant of its centralizer
elementary equivalent to $PSL_2(K)\times PSL_2(K)$.
\end{proof}

Now we need only to consider groups of the type $B_2$, i.\,e. groups $EO_5(K)$, 
that contain only the involutions conjugate to $
daig[-1,-1,-1,-1,1]$, with the central quotient of the commutant 
of its centralizer isomorphic to $PEO_4(K)$. This group can be
elementary equivalent to $PSL_2(K)\times
PSL_2(K)$, so we need to find some sentences separating groups
of the types $A_3$ and  $B_2$.

\begin{lemma}\label{klass4}
There exists a first order sentence $\varphi_{A_3-B_2}$, that holds in any
adjoint groups of the type $A_3$ and does not hold in all 
adjoint Chevalley groups of the type~$B_2$. 
\end{lemma}
\begin{proof}
It is clear that in one of the groups $G(A_3)$ and $G(B_2)$ there exists
an involution that is not conjugate to  the involution mentioned above,
then we can separate these groups. So we need to consider the case
when in these groups there is only one conjugate class of involutions.

Let us consider a matrix $M$, satisfying the following formula:
\begin{multline*}
DDiag_2(M):=\exists I (I^2=1\land I\ne 1\land \exists M_1\exists M_2
(M_1 I=IM_1\land\\
\land M_2I=IM_2\land M=M_1M_2M_1^{-1}M_2^{-1}\land\\
\land \forall X\forall Y (XM=MY\land YM=MY\Rightarrow XY=YX)\land\\
\land \forall N\forall X (NI=IN\land 
NM=MN\land (XI=IX\land M=XNX^{-1})\Rightarrow\\
\Rightarrow X^2M=MX^2)).
\end{multline*}

This formula states that there exists such an involution $I$ that

1) $M$ lies in the commutant of its centralizer;

2) any matrices commuting with~$M$, commutes;

3) for any matrix~$N$, commuting with $M$ and~$I$ and conjugate to~$M$
by the matrix~$X$, that commutes with~$I$, the condition
$X^2M=MX^2$ holds.

Let us look what matrices satisfy the formula $DDiag_2(M)$ in the group
 $G(A_3)=PSL_4(K)$.

An involution $I$ has in some basis the form
$diag[\xi,\xi,-\xi,-\xi]$, $\xi^4=1$.
Therefore in this basis the matrix  $M$ has the form
$$
M=\begin{pmatrix}
a& b& 0& 0\\
c& d& 0& 0\\
0& 0& f& e\\
0& 0& g& h
\end{pmatrix}.
$$
The second part of the formula in this case means only that the matrix
 $M$ has no Jordan cells with the same proper values.
Thus, every cell
$$
\begin{pmatrix}
a& b\\
c& d
\end{pmatrix}\text{ and } 
\begin{pmatrix}
f& e\\
g& h
\end{pmatrix}
$$
either is diagonalizable in~$K^*$ with different proper values,
or is a Jordan cell $2\times 2$. 

Let 
$$
\begin{pmatrix}
a& b\\
c& d
\end{pmatrix}
$$
be a Jordan cell
$$
\begin{pmatrix}
\alpha& 1\\
0& \alpha
\end{pmatrix}.
$$
Choose
$$
N=\begin{pmatrix}
\alpha& 2& 0& 0\\
0& \alpha & 0& 0\\
0& 0& 1& 0\\
0& 0& 0& 1
\end{pmatrix},\quad
X= \begin{pmatrix}
2& 0& 0& 0\\
0& 1& 0& 0\\
0& 0& 1/2& 0\\
0& 0& 0& 1
\end{pmatrix}
.$$
These matrices contradict to the third part of the formula $DDiag_2(M)$.

Therefore, both cells are diagonalizable, and all proper values
a different. Thus, in the group $PSL_4(K)$ the formula $DDiag_2(M)$ defines
all diagonalizable in~$\overline K$ matrices with different
diagonal elements, and only them.

Now let us look what matrices satisfy  this formula in the group
$G(B_2)=EO_5(K)$.
An involution $I$ has the form $diag[-1,-1,-1,-1,1]$ in a basis, where
a form is
$$
\begin{pmatrix}
Q_1& 0\\
0& q_2
\end{pmatrix},\quad Q_1\in GL_4(K), \ Q_1^T=Q_1, \ q_2\in K^*.
$$
Therefore, 
$$
M=\begin{pmatrix}
M_1& 0\\
0& 1
\end{pmatrix},\quad M_1\in SO_4(K).
$$
Let us take its Jordan form:
$$
\widetilde M=\begin{pmatrix}
\widetilde{M_1}&0\\
0& 1
\end{pmatrix},\quad
\widetilde Q=\begin{pmatrix}
\widetilde{Q_1}& 0\\
0& 1
\end{pmatrix}.
$$

Consider possible variants:

1) if $\widetilde M$ is diagonal  with different diagonal elements, 
then it is easy to show that it satisfies the formula $DDiag_2(M)$.

2) The matrix $\widetilde M=diag[\alpha,\alpha,\beta  ,\gamma,\delta]$
in any case can not satisfy the formula $DDiag_2(M)$.

3) If 
$$
\widetilde M=\begin{pmatrix}
\alpha& 1& 0& 0& 0\\
0& \alpha& 0& 0& 0\\
0& 0& \beta& 0& 0\\
0 & 0& 0& \gamma& 0\\
0& 0& 0& 0& 1
\end{pmatrix},
$$
then we have the conditions
$$
\begin{pmatrix}
\alpha& 1& 0& 0 \\
0& \alpha& 0& 0\\
0& 0& \beta& 0\\
0 & 0& 0& \gamma
\end{pmatrix}
\begin{pmatrix}
q_1& q_2& q_3& q_4 \\
q_2& q_5& q_6& q_7\\
q_3& q_6& q_8& q_9\\
q_4 & q_7& q_9& q_{10}
\end{pmatrix}
\begin{pmatrix}
\alpha& 0& 0& 0 \\
1& \alpha& 0& 0\\
0& 0& \beta& 0\\
0 & 0& 0& \gamma
\end{pmatrix}=
\begin{pmatrix}
q_1& q_2& q_3& q_4 \\
q_2& q_5& q_6& q_7\\
q_3& q_6& q_8& q_9\\
q_4 & q_7& q_9& q_{10}
\end{pmatrix}.
$$
So we have the system of equations
$$
\begin{cases}
\alpha^2 q_1+\alpha q_2+\alpha q_2+q_5=q_1,\\
\alpha^2 q_2+\alpha q_5=q_2,\\
\alpha \beta q_3+\beta q_6=q_3\\
\alpha \gamma q_4+\gamma q_7=q_4,\\
\alpha^2 q_5=q_5,\\
\alpha \beta q_6=q_6,\\
\alpha \gamma q_7=q_7,\\
\beta^2 q_8=q_8,\\
\gamma \beta q_9=q_9,  \\
\gamma^2 q_{10}=q_{10}.
\end{cases}
$$

Now we also need to consider different possible cases
for different $\alpha,\beta,\gamma$.

a) $\alpha =\pm 1$, therefore $q_5=0\Rightarrow q_2=0\Rightarrow (q_6\ne 0\lor
q_7\ne 0)\Rightarrow \beta=\alpha\lor \gamma =\alpha$. 
If $\beta=\alpha =\pm 1$, $\gamma \ne \alpha$, then $q_6=q_7=0$, and it is impossible.
If $\gamma=\beta=\pm 1$, then the matrix does not satisfy the formula 
$DDiag_2(M)$.

b) if $\alpha\ne \pm 1$, then $q_1=q_2=q_5=0$. Therefore either $q_6\ne 0$, or
$q_7\ne 0$. It means that either $\beta=1/\alpha$, or $\gamma=1/\alpha$.
Let $\beta=1/\alpha$, $\gamma\ne 1/\alpha$. Then $q_7=0$, $q_6=0$,
and it is impossible.

Consequently, the matrix $\widetilde M$ can not have this form.

4) Let
$$
\widetilde M=\begin{pmatrix}
\alpha& 1& 0& 0\\
0& \alpha& 0& 0\\
0& 0& \beta& 1\\
0& 0& 0& \beta
\end{pmatrix}.
$$
From the conditions
$$
\begin{pmatrix}
\alpha& 1& 0& 0 \\
0& \alpha& 0& 0\\
0& 0& \beta& 1\\
0 & 0& 0& \beta
\end{pmatrix}
\begin{pmatrix}
q_1& q_2& q_3& q_4 \\
q_2& q_5& q_6& q_7\\
q_3& q_6& q_8& q_9\\
q_4 & q_7& q_9& q_{10}
\end{pmatrix}
\begin{pmatrix}
\alpha& 0& 0& 0 \\
1& \alpha& 0& 0\\
0& 0& \beta& 0\\
0 & 0& 1& \beta
\end{pmatrix}=
\begin{pmatrix}
q_1& q_2& q_3& q_4 \\
q_2& q_5& q_6& q_7\\
q_3& q_6& q_8& q_9\\
q_4 & q_7& q_9& q_{10}
\end{pmatrix}
$$
we have the system of equations
$$
\begin{cases}
\alpha^2 q_1+\alpha q_2+\alpha q_2+q_5=q_1,\\
\alpha^2 q_2+\alpha q_5=q_2,\\
\alpha \beta q_3+\beta q_6+\alpha q_4+q_7=q_3\\
\alpha \beta q_4+\beta q_7=q_4,\\
\alpha^2 q_5=q_5,\\
\alpha \beta q_6+\alpha q_7=q_6,\\
\alpha \beta q_7=q_7,\\
\beta^2 q_8+2\beta q_9+q_{10}=q_8,\\
\beta^2 q_9+\beta q_{10}=q_9,  \\
\beta^2 q_{10}=q_{10}.
\end{cases}
$$

a) $\alpha,\beta\ne \pm 1$, $\alpha\beta\ne 1$ $\Rightarrow$
$q_7=q_{10}=0$, therefore $q_2=0\Rightarrow q_6=0$, but it is impossible.

b) $\alpha,\beta\ne \pm 1$, $\alpha \beta=1$.
Then $q_5=q_{10}=0\Rightarrow q_9=q_2=q_7=q_1=q_8=0$, therefore
$$
\widetilde Q=\begin{pmatrix}
0& 0& q& \beta\\
0& 0& \alpha& 0\\
q& \alpha& 0& 0\\
\beta& 0& 0& 0
\end{pmatrix}.
$$
Matrices in this basis satisfy the condition
$$
\begin{pmatrix}
A& B\\
C& D
\end{pmatrix}
\widetilde Q \begin{pmatrix}
A^T& C^T\\
B^T& D^T
\end{pmatrix}=Q.
$$
Thus,
in this basis we have matrices
$$
\begin{pmatrix}
A& 0\\
0& D
\end{pmatrix},\quad D=\begin{pmatrix}
q& \alpha\\
\beta& 0
\end{pmatrix}
(A^T)^{-1}
\begin{pmatrix}
q& \alpha\\
\beta& 0
\end{pmatrix}^{-1}, 
$$
$A$ is an arbitrary matrix.

Consequently, the matrix $\widetilde M$ does not satisfy the last part
of the formula $DDiag_2(M)$.

c) $\alpha=\pm 1$, $\beta\ne \pm 1$. Then $q_7=q_{10}=q_9=q_4=0$, 
but it is impossible.

d) $\alpha,\beta=\pm 1$, $\alpha\beta=-1$. Then
$q_7=q_4=q_6=q_3=q_5=q_2=0$,
it is impossible.

e) $\alpha,\beta=\pm 1$, $\alpha \beta=1$. Then we have the following system
of equations:
$$
\begin{cases}
 q_1\pm 2 q_2+q_5=q_1,\\
 q_2\pm  q_5=q_2,\\
 q_3\pm q_6\pm q_4+q_7=q_3\\
 q_4\pm q_7=q_4,\\
 q_5=q_5,\\
q_6\pm q_7=q_6,\\
 q_7=q_7,\\
 q_8\pm 2 q_9+q_{10}=q_8,\\
q_9\pm q_{10}=q_9,  \\
 q_{10}=q_{10}.
\end{cases}
$$

Therefore, $q_5=q_7=q_{10}=q_2=q_9=0$, $q_6=-q_4$, i.\,e.
$$
\widetilde Q=\begin{pmatrix}
q_1& 0& q_3& q_4\\
0& 0& -q_4& 0\\
q_3& -q_4& q_8& 0\\
q_4& 0& 0& 0
\end{pmatrix}.
$$
We can directly show that the matrices
$$
\begin{pmatrix}
0& b& 1& 0\\
0& 0& 0& 1\\
-1& f& 0& h\\
0& -1& 0& 0
\end{pmatrix},\quad b=\frac{q_8-q_1}{2q_4},\ f=\frac{-2q_3}{q_4},\
h=\frac{q_1-q_8}{2q_4}
$$
and
$$
\begin{pmatrix}
0& b'& -1/2& d\\
0& 0& 0& -1/2\\
2& 0& 0& h'\\
0& 2& 0& 0
\end{pmatrix},\quad
b'=\frac{q_1-q_8/4}{q_4},\ h'=\frac{q_8-4q_1}{4q_4},\ d=\frac{q_4}{q_3}
$$
belong to our group, do not commute, but each of them commutes
with $\widetilde M$. Therefore, the matrix $\widetilde M$ does not satisfy
the initial formula.

Therefore, we have shown that in  both groups the formula $DDiag_2(M)$ 
defines diagonalizable in~$\overline K$ matrices with different
elements on the diagonal.
Matrices that commute with~$M$, form the subgroup~$H$
of diagonalizable in the same basis matrices, 
and the quotient group of the normalizer of~$H$ by~$H$
is isomorphic to the Weil group~$W$ of~$G$.
In the first case it is the group $S_4$, in the second case it is
the product of  $S_2$ and $(\mathbb Z_2)^2$. Therefore it is possible to write 
a sentence, that separates the groups of types $A_3$ and $B_2$.
\end{proof}

Therefore, according to Lemmas~\ref{klass3} and~\ref{klass4},
we found the sentences, separating adjoint groups of the type $A_3$
from other classical groups.

\begin{lemma}\label{klass5}
For every $l\ge 4$ there exists a first order sentence $\varphi_{A_l}$
that holds in any adjoint groups of the type $A_l$ and does not hold
in all adjoint Chevalley groups of other classical types.
\end{lemma}
\begin{proof}
All Chevalley groups $A_l$, $l\ge 4$  are characterized by the fact
that they have at least two non-conjugate involutions
with infinite centers of centralizers (with center of centralizers of order~$>2$).

Then we can easily separate them from each other with the help of
the groups of smaller dimensions 
(the central quotient of the commutant of the centralizer of an involution
is either $PSL_k(K)\times PSL_{m-k}(K)$, or
$PSL_{m-1}(K)$).
\end{proof}

\begin{lemma}\label{klass6}
For every  $l\ge 2$ there exists a first order sentence $\varphi_{B_l}$
that holds in any
adjoint group of the type $B_l$ and does not hold
in all adjoint Chevalley groups of other classical types. 
\end{lemma}
\begin{proof}
The group $EO_5(K)$ (i.\,e., an adjoint group of the type $B_2$)
is already separated from the other groups
(see Lemma~\ref{klass4}).

Let us consider the group $EO_7(K)$ (the type $B_3$).
If this group contains an involution $diag[-1,-1,1,\dots,1]$, 
then it is uniquely characterized by the fact that it contains
an involution with the central quotient 
of the commutant of its centralizer
 $EO_5(K)$.

If there is no such an involution, then this group contains only one 
conjugate class: $diag[-1,-1,-1,-1,1,1,1]$. 
We know that the given group is not of the type $A_l$ (according to 
Lemmas \ref{klass1}--\ref{klass5}),
and not of the type $B_2$. It can not have the type $B_l$ ($l>3$), since
in this case there are more than one conjugate classes of involutions. 
It can not have the type $C_l$ ($l\ge 3$) or $D_l$ ($l\ge 4$). Therefore
it is the group of the type $B_3$.

Now let us consider  $l> 3$.

In such groups there are involutions
(of the form $diag[-1,-1,-1,-1,1,\dots,1]$), with central quotient
of the commutant of its centralizer  $EO_{2l-3}(K)\times EO_4(K)$.
We already have the sentence characterizing $EO_{2l-3}(K)$. 
Since this group never appears in other cases, then we have separate the
type $B_l$.
\end{proof}

\begin{lemma}\label{klass7}
For every $l\ge 3$ there exists a first order sentence $\varphi_{C_l}$
that holds in any 
adjoint group of the type $C_l$ and does not hold in all
adjoint Chevalley groups of other classical types.
\end{lemma}
\begin{proof}
At first we consider the groups $PSp_6(K)$ (i.\,e., the groups of the type
 $C_3$).

According to Lemmas \ref{klass1}--\ref{klass6} we only need to separate
a group of the type $C_3$ from  groups of the types  $C_l$ and $D_l$, $l\ge 4$. 
It contains an involution with central quotient of the commutant of it
centralizer isomorphic to $PSL_3(K)$.
Other groups can not satisfy this property.

Now let us consider the general case of the groups $PSp_{2l}(K)$ ($l\ge 4$).

We can act recursively:
there exists an involution with the central quotient
of the commutant of its centralizer isomorphic to
$$
PSp_{2(l-1)}(K)\times PSL_2(K).
$$
\end{proof}

\begin{lemma}\label{klass8}
For every $l\ge 4$ there exists a first order sentence $\varphi_{D_l}$
that holds in any
adjoint group of the type $D_l$ and does not hold in all adjoint
Chevalley groups of other classical types. 
\end{lemma}
\begin{proof}
According to Lemmas \ref{klass1}--\ref{klass7}
we only need to separate group of the types $D_l$ and $D_m$ for
$l\ne m$.

For even~$l$
there exists an involution with the central quotient
of the commutant of its centralizer isomorphic to $PSL_l(K)$, 
and for odd~$l$ it is
$PEO_{2(l-1)}(K)$.

Therefore we can separate all group of the type $D_l$.
\end{proof}

From Lemmas \ref{klass1}--\ref{klass8} it follows that
no two adjoint elementary Chevalley groups with different 
root systems can be elementary equivalent.

In the next sections we shall state this result for
exceptional root systems.

\section{Chevalley groups of the type $G_2$.}\leavevmode

To study a Chevalley groups of the type $G_2$, we need to consider
the Lie algebra of the type $G_2$. The roots of this algebra can be considered
as linear combinations of basis vectors of $3$-dimensional vector space.
Namely,
$\alpha_1=e_1-e_2$, $\alpha_2=-e_1+e_2+e_3$, $\alpha_3=\alpha_1+\alpha_2=
e_3-e_1$, $\alpha_4=2\alpha_1+\alpha_2=e_3-e_2$,  $\alpha_5=3\alpha_1+\alpha_2=
-2e_2+e_1+e_3$, $\alpha_6=3\alpha_1+2\alpha_2=2e_3-e_1-e_2$. Therefore
we have two simple roots $\alpha_1$ and $\alpha_2$, and six positive
roots $\alpha_1,
\dots,\alpha_6$, the general number of roots is~$12$.  
Now let us write a table
with all values $\langle \alpha,\beta\rangle =2(\alpha,\beta)/(\beta,\beta)$
for simple roots $\beta$ and positive roots~$\alpha$.

{\small
\begin{center}

\begin{tabular}{|c|c|c|c|c|c|}
\hline
 & $\alpha_1$ & $\alpha_2$&
 & $\alpha_1$ & $\alpha_2$\\
\hline
$\alpha_1$ & $2$ & $-1$&
$\alpha_4$ & $1$ & $0$\\
\hline
$\alpha_2$ & $-3$ & $2$&
$\alpha_5$ & $3$ & $-1$\\
\hline
$\alpha_3$ & $-1$ & $1$&
$\alpha_6$ & $0$ & $1$\\
\hline
\end{tabular}

\end{center}
}

There exists a faithful representation of the Lie algebra $G_2$ in
$SL_7(\mathbb C)$. To show it we only need to write the images
of the elements $X_\alpha$ for the roots
$\alpha$ of $G_2$. We denote  by $X_i$ 
the matrix corresponding to the element $X_{\alpha_i}$ for a positive
root~$\alpha_i$, and by $X_{-i}$ the matrix corresponding
to a negative root $-\alpha_i$. Here are these  12 matrices:
$X_1=-ie_{13}+2ie_{21}+ie_{46}-ie_{75}$,
$X_{-1}=-ie_{12}+2ie_{31}-ie_{57}+ie_{64}$,
$X_2=e_{37}-e_{62}$,
$X_{-2}=-e_{26}+e_{73}$,
$X_3=ie_{17}-ie_{35}+ie_{42}-2ie_{61}$,
$X_{-3}=-ie_{16}-ie_{24}+ie_{53}+2ie_{71}$,
$X_{4}=-ie_{15}+i2_{27}+2ie_{41}-ie_{63}$,
$X_{-4}=-ie_{14}-ie_{36}+2ie_{51}+ie_{72}$,
$X_5=-e_{27}+e_{43}$,
$X_{-5}=e_{34}-e_{72}$,
$X_6=-e_{47}+e_{65}$,
$X_{-6}=e_{56}-e_{74}$.

By this Lie algebra we can construct the Chevalley group of the type $G_2$.
In the case $i\in K$ there exists a representation of these group
in the group $GL_7(K)$, in the other case it is in the group $GL_{14}(K)$.

With the help of $x_1(t),\dots, x_6(t), x_{-1}(t),\dots, x_{-6}(t)$ 
we find the matrices corresponding to $w_1(t),\dots, w_6(t)$ and
$h_1(t),\dots,h_6(t)$:

$w_1(t)=-e_{11}+t^2e_{23}+t^{-2}e_{32}+ite_{46}-it^{-1}e_{57}+it^{-1}
e_{64}-ite_{75}$,
$h_1(t)=[1,t^2,1/t^2,t,1/t,1/t,t],$
$w_2(t)=e_{11}+t^{-1}e_{26}+t3_{37}+e_{44}+e_{55}-te_{62}-t^{-1}e_{73}$,
$h_2(t)=[1,1/t,t,1,1,t,1/t]$,
$w_3(t)=-e_{11}-it^{-1}e_{24}+ite_{35}-ite_{42}+it^{-1}e_{53}+t^2 e_{67}
+t^{-2} e_{76}$,
$h_3(t)=[1,1/t,t,1/t,t,t^2,1/t^2],$
$w_4(t)=-e_{11}-ite_{27}+it^{-1}e_{36}+t^2 e_{45} +t^{-2} e_{54} +it e_{63}
-it^{-1}e_{72}$,
$h_4(t)=[1,t,1/t,t^2,1/t^2,t,1/t],$
$w_5(t)=e_{11}-te_{25}-t^{-1}e_{34}+te_{43}+t^{-1}e_{52}+e_{66}+e_{77}$,
$h_5(t)=[1,t,1/t,t,1/t,1,1],$
$w_6(t)=e_{11}+e_{22}+e_{33}-te_{47}-t^{-1}e_{56}+te_{65}+t^{-1}e_{74}$,
$h_6(t)=[1,1,1,t,1/t,t,1/t].$

Therefore, the subgroup $H$ of $G$ (that coincides with the torus~$T$,
since $G$ is universal) consists of diagonal matrices
with proper values $[1, u, u^{-1}, uv, u^{-1}v^{-1},
v, v^{-1}]$. In this group there are only three involutions, they have the form
 $[1,1,1,-1,-1,-1,-1]$, $[1, -1,-1,-1,-1,1,1]$ 
and $[1,-1,-1,1,1,-1,-1]$. 

The Weil group of this Chevalley group is isomorphic
to the dihedral group of the order~$12$.
It is represented by the matrices
$e$, $w_1(1)$, $w_2(1)$, $w_3(1)$,
$w_4(1)$, $w_5(1)$, $w_6(1)$,
$w_1(1)w_2(1)=-e_{11}+e_{27}-e_{36}-ie_{42}+ie_{53}+ie_{64}-ie_{75}$,
$w_1(1)w_3(1)=e_{11}+ie_{25}-ie_{34}+ie_{47}-ie_{56}+e_{62}+e_{73}$,
$w_1(1)w_4(1)=e_{11}+ie_{26}-ie_{37}-e_{43}-e_{52}+ie_{65}-ie_{74}$,
$w_1(1)w_5(1)=-e_{11}-e_{24}-e_{35}-ie_{46}+ie_{57}-ie_{63}+ie_{72}$,
$w_1(1)w_6(1)=-e_{11}+e_{23}+e_{32}-ie_{45}+ie_{54}+ie_{67}-ie_{76}$.

\begin{lemma}\label{razlG}
There exists a first order sentence $\varphi_{G_2}$
that holds in any adjoint group of the type
 $G_2$ and does not hold in all  classical 
adjoint Chevalley groups. 
\end{lemma}

\begin{proof}
From the representation of the Weil group we see that
its element $w_1(1)w_6(1)$ is an involution and commutes
with all involutions of the group~$H$, 
and such an element of the Weil group is unique.
It is  clear that its products with involutions of the group~$H$
also commutes with all Weil group elements and these products are involutions.

Let us show that the elements $h_2(-1)$ and $w_1(1)w_6(1)h_2(-1)$ 
are conjugate.

Actually, 
\begin{multline*}
x_{-6}(-1/2)h_2(-1)x_{-6}(1/2)=x_{-6}(-1)h_2(-1),\\
x_6(1)x_{-6}(-1)h_2(-1)x_6(-1)=
(x_6(1)x_{-6}(-1)x_6(-1))(x_6(1)h_2(-1)x_6(-1))=\\
x_6(1)x_{-6}(-1) x_6(-1)x_6(2)h_2(-1)=\\
=x_6(1)x_{-6}(-1)x_6(1)h_2(-1)=w_6(1)h_2(-1),\\
 x_{-1}(-1/2)w_6(1)h_2(-1)
x_{-1}(1/2)=x_{-1}(-1)w_6(1)h_2(-1),\\
x_1(1)x_{-1}(-1)w_6(1)h_2(-1)x_1(-1)=w_1(1)w_6(1)h_2(-1).
\end{multline*}

Therefore, it is easy to see that the elements
 $h_1(-1)$, $h_2(-1)$, $h_3(-1)$,
$w_1(1)w_6(1)h_1(-1)$, $w_1(1)w_6(1)h_2(-1)$ and $w_1(1)w_6(1)h_3(-1)$ 
are conjugate.

Besides it, if $i\in K$, 
then the elements $w_1(1)w_6(1)$ and $w_1(1)w_6(1)h_1(-1)$
are conjugate by $h_1(i)$, i.\,e.
all mentioned involutions.

Therefore, we have a set of $7$ commuting involutions, where either
all seven, or six involutions are conjugate. It is clear that all
centralizers of these involutions are isomorphic.
Let us look what is the centralizer of $h_1(-1)$.

Let $x\in G$ belong to the centralizer of $h_1=h_1(-1)$.

From the Bruhat decomposition (see~\cite{Steinberg})
it follows that any element of~$G$
can be uniquely represented in the form
$$
x_1(t_1)x_2(t_2)x_3(t_3)x_4(t_4)x_5(t_5)x_6(t_6)hw x_1(u_1)x_2(u_2)x_3(u_3)
x_4(u_4)x_5(u_5)x_6(u_6),
$$
where $h\in H$, $w$ is one of the elements of the Weil group,
$u_i=0$ for all~$i$ with $w(\alpha_i)\in \Phi^+$.
Since $h_1 xh_1^{-1}=x$, we have that
\begin{multline*}
x_1(t_1)x_2(-t_2)x_3(-t_3)x_4(-t_4)x_5(-t_5)x_6(t_6)h(h_1(-1)wh_1(-1))\times\\
\times x_1(u_1)x_2(-u_2)x_3(-u_3)
x_4(-u_4)x_5(-u_5)x_6(u_6)=\\
=x_1(t_1)x_2(-t_2)x_3(-t_3)x_4(-t_4)x_5(-t_5)x_6(t_6)hh'w\times\\
\times  x_1(u_1)x_2(-u_2)x_3(-u_3)
x_4(-u_4)x_5(-u_5)x_6(u_6)=\\
=x_1(t_1)x_2(t_2)x_3(t_3)x_4(t_4)x_5(t_5)x_6(t_6)hw\times\\
\times x_1(u_1)x_2(u_2)x_3(u_3)
x_4(u_4)x_5(u_5)x_6(u_6),
\end{multline*}
where $h'$ is some element of the group $H$. 
The uniqueness of the Bruhat decomposition implies
 $t_2=t_3=t_4=t_5=u_2=u_3=u_4=u_5=0$, $h'=1$, i.\,e. 
$h_1(-1)wh_1(-1)=w$. It is clear to check that either
$w=1$, or $w=w_1(1)$, or $w=w_6(1)$, or $w=w_1(1)w_6(1)$.
Therefore the centralizer of $h_1(-1)$ is the product
of the group $H$, the group $X_1$, generated by $x_1(t)$ and $x_{-1}(t)$,
and the group $X_6$, generated by $x_6(t)$ and $x_{-6}(t)$.
The commutant of the considered group is generated by the groups
 $X_1$ and $X_6$.
The central quotient of this commutant is
 $PSL_2(K)\times PSL_2(K)$. 

As we remember, besides the group $G_2$, only Chevalley groups 
 $PSL_4(K)$ and $EO_5(K)$ can satisfy this property.

We have already shown that the commutant of the centralizer of $h_1(-1)$
is generated by the groups $X_1$ and $X_6$. Similarly, 
the commutant of the centralizer of $h_2(-1)$ is generated by the groups
 $X_2$ and $X_4$, 
and the commutant of the centralizer of $h_3(-1)=h_1(-1)h_2(-1)$
is generated by $X_3$  and $X_5$.  Note that the pairs of groups
$X_1$ and $X_2$, $X_1$ and $X_3$, $X_1$ and $X_4$, $X_1$ and $X_5$,
$X_2$ and $X_3$, $X_3$ and $X_4$, $X_3$ and $X_6$, $X_4$ and $X_5$,
$X_4$ and $X_6$ generate the whole group~$G$, and the pairs $X_2$ and $X_5$,
$X_2$ and $X_6$, $X_5$ and $X_6$ generate only a proper subgroup
of~$G$, that is a Chevaley group of the type $A_2$. 

Neither the group $PSL_4(K)$, nor the group $EO_5(K)$ can satisfy this property.
It is clear that this property can be written  as a first order
sentence, therefore
the Chevalley group $G_2$ can not be elementary equivalent to any
classical Chevalley group.
\end{proof}

\section{Chevalley groups of the type $F_4$.}\leavevmode

In the root system $F_4$ there are $48$ roots,  $24$ roots are positive.
These positive roots are generated by $4$ simple roots
$\alpha_1=e_2-e_3$, $\alpha_2=e_3-e_4$, $\alpha_3=e_4$ and
$\alpha_4=\frac{1}{2}(e_1-e_2-e_3-e_4)$. The first two roots are long,
the second two are short. These are other positive roots:
{\small
$\alpha_5=\alpha_1+\alpha_2$, $\alpha_6=\alpha_2+\alpha_3$, $\alpha_7=
\alpha_3+\alpha_4$, $\alpha_8=\alpha_1+\alpha_2+\alpha_3$, $\alpha_9=
\alpha_2+\alpha_3+\alpha_3$, $\alpha_{10}=\alpha_2+\alpha_3+\alpha_4$, 
$\alpha_{11}=\alpha_1+\alpha_2+\alpha_3+\alpha_4$, $\alpha_{12}=
\alpha_2+\alpha_3+\alpha_3+\alpha_4$, $\alpha_{13}=\alpha_1+\alpha_2+\alpha_3+\alpha_3$,
$\alpha_{14}=\alpha_1+\alpha_2+\alpha_2+\alpha_3+\alpha_3$,
$\alpha_{15}=\alpha_1+\alpha_2+\alpha_3+\alpha_3+\alpha_4$,
$\alpha_{16}=\alpha_2+\alpha_3+\alpha_3+\alpha_4+\alpha_4$,
$\alpha_{17}=\alpha_1+\alpha_2+\alpha_2+\alpha_3+\alpha_3+\alpha_4$,
$\alpha_{18}=\alpha_1+\alpha_2+\alpha_2+\alpha_3+\alpha_3+\alpha_3+\alpha_4$,
$\alpha_{19}=\alpha_1+\alpha_2+\alpha_2+\alpha_3+\alpha_4+\alpha_4$,
$\alpha_{20}=\alpha_1+\alpha_2+\alpha_2+\alpha_3+\alpha_3+\alpha_3+
\alpha_4+\alpha_4$, $\alpha_{21}=\alpha_1+\alpha_1+\alpha_2+\alpha_2+\alpha+2+
\alpha_3+\alpha_3+\alpha_3+\alpha_4+\alpha_4$, $\alpha_{22}=\alpha_1+\alpha_2+\alpha_3+
\alpha_3+\alpha_4+\alpha_4$, $\alpha_{23}=\alpha_1+\alpha_2+\alpha_2+\alpha_2+
\alpha_3+\alpha_3+\alpha_3+\alpha_3+\alpha_4+\alpha_4$, $\alpha_{24}=
\alpha_1+\alpha_2+\alpha_2+\alpha_3+\alpha_3+\alpha_3+\alpha_3+\alpha_4+\alpha_4$.
}

\begin{lemma}\label{razlF}
There exists a first order sentence $\varphi_{F_4}$
that holds in any adjoint 
groups of the type $F_4$ and does not hold in all classical
adjoint Chevalley groups, and also in Chevalley groups of the type~$G_2$. 
\end{lemma}

\begin{proof}
Let us find the centralizer of $h_1(-1)$.

Let $Z(h_1)\ni x=x_1(t_1)\dots x_{24}(t_{24})hw x_1(u_1)\dots x_{24}(u_{24})$,
where $h\in H$, $w$ is a product of some $w_i(1)$.
Then 
$$
h_1(-1) x h_1(-1)=
x_1(\pm t_1)\dots x_{24}(\pm t_{24}) h h' w
x_1(\pm u_1)\dots x_{24}(\pm u_{24}),
$$
 and $h_1(-1) x_i(s) h_1(-1)=x_i(s)$
if and only if $i=1,3,4,7,14,17,18,19,20,24$. Therefore,
\begin{multline*}
x=x_{\alpha_1}(t_1)x_{\alpha_3}(t_3)x_{\alpha_4}(t_4)
 x_{\alpha_3+\alpha_4}(t_7)
x_{\alpha_{14}}(t_{14})\times\\
x_{\alpha_{14}+\alpha_4}(t_{17})x_{\alpha_{14}+\alpha_3+
\alpha_4}(t_{18})x_{\alpha_{14}+x_{\alpha_4}+x_{\alpha_4}}(t_{19})
x_{\alpha_{14}+
\alpha_3+\alpha_4+\alpha_4}(t_{20})\times\\
\times x_{\alpha_{14}+\alpha_3+\alpha_3+
\alpha_4+\alpha_4}(t_{24}) h w x_1(u_1)\times\\
\times x_3(u_3)x_4(u_4)x_7(u_7)x_{14}(u_{14})
x_{17}(u_{17})x_{18}(u_{18})x_{19}(u_{19})x_{20}(u_{20})x_{24}(u_{24}).
\end{multline*}
Therefore, $h_1(-1)$ has to commute with $w$. We know that
$w_{\alpha}(1) h_\beta(-1) w_\alpha(1)^{-1}=h_{w_\alpha(\beta)}(-1)$.
Consequently, $w_{i_1}(1)\dots w_{i_m}(1)$ commutes with  $h_1(-1)$
if and only if $w_{\alpha_{i_1}} \dots w_{\alpha_{i_m}}(\alpha_1)=
\pm \alpha_1$. 

It is clear that any product of $w_1(1)$, $w_3(1)$, $w_4(1)$
и $w_{14}(1)$ commutes with $h_1(-1)$. 
Easy to check that other elements of the Weil group
do not commute with $h_1(-1)$. Thus, the centralizer of 
$h_1(-1)$ is a product of the group $H$ and the subgroup of~$G$, generated
by $X_1$, $X_3$, $X_4$ and $X_{14}$,
i.\,e., it is a product of $H$,  the Chevalley groups of the type
$A_1$ and the Chevalley groups of the type $C_3$. As above we can show that
its commutant 
is a central product of the Chevalley of the type $A_1$
and the Chevalley groups of the type $C_3$, i.\,e., the groups
 $SL_2(K)$ and $Sp_6(K)$,
and the central quotient of this commutant is
$PSL_2(K)\times PSp_6(K)$. 

Similarly, the centralizer of $h_3(-1)$
is generated by $H$ and the subgroups $X_1$, $X_2$, $X_3$, $X_{16}$,
i.\,e., is a product of $H$ and the Chevalley group of the type
$B_4$, and the central quotient of this commutant is $EO_9(K)$.

To write a sentence separating  Chevalley groups of the type $F_4$
from all considered Chevalley groups, we only can mention that
we have no other Chevalley groups
with a central quotient of a commutant of a centralizer of some
involution isomorphic to $PSL_2(K)\times PSp_6(K)$, and with
a central quotient of a commutant of  a centralizer of some other
involution isomorphic to $EO_9(K)$.
\end{proof}

\section{Chevalley groups of the type $E_6$.}\leavevmode

In the root system $E_6$ there are $72$ roots,  $36$ roots are positive.
These positive roots are generated by six simple roots
$\alpha_1=\frac{1}{2}(e_1+e_8-e_2-e_3-e_4-e_5-e_6-e_7)$, $\alpha_2=e_1+e_2$, 
$\alpha_3=e_2-e_1$, $\alpha_4=e_3-e_2$, $\alpha_5=e_4-e_3$ and
$\alpha_6=e_5-e_4$.  These are other positive roots
(for simplicity we write $i$ instead of $\alpha_i$):
{\small $\alpha_7=1+3$, $\alpha_8=3+4$, $\alpha_9=
4+5$, $\alpha_{10}=5+6$, 
$\alpha_{11}=1+3+4$, $\alpha_{12}=
3+4+5$, $\alpha_{13}=2+4$,
$\alpha_{14}=2+3+4$,
$\alpha_{15}=2+4+5$,
$\alpha_{16}=1+2+3+4$,
$\alpha_{17}=2+3+4+5$,
$\alpha_{18}=1+3+4+5$,
$\alpha_{19}=3+4+5+6$,
$\alpha_{20}=2+4+5+6$, 
$\alpha_{21}=2+3+4+4+5$, 
$\alpha_{22}=2+3+4+5+6$, 
$\alpha_{23}=1+2+3+4+5$, 
$\alpha_{24}=1+3+4+5+6$,
$\alpha_{25}=1+2+3+4+4+5$,
$\alpha_{26}=1+2+3+3+4+4+5$,
$\alpha_{27}=2+3+4+4+5+6$,
$\alpha_{28}=2+3+4+4+5+5+6$,
$\alpha_{29}=1+2+3+3+4+4+5+6$,
$\alpha_{30}=1+2+3+4+5+6$,
$\alpha_{31}=1+2+3+4+4+5+6$,
$\alpha_{32}=1+2+3+4+4+5+5+6$,
$\alpha_{33}=4+5+6$,
$\alpha_{34}=1+2+3+3+4+4+5+5+6$,
$\alpha_{35}=1+2+3+3+4+4+4+5+5+6$,
$\alpha_{36}=1+2+2+3+3+4+4+4+5+5+6$.}

\begin{lemma}\label{razlE6}
There exists a first order sentence $\varphi_{E_6}$
that holds in any
adjoint group of the type $E_6$ and does not hold
in all classical adjoint Chevalley groups and also in all Chevalley groups
of types $G_2$ and~$F_4$. 
\end{lemma}

\begin{proof}
Let us find the centralizer of $h_1$.

Let $Z(h_1)\ni x=x_1(t_1)\dots x_{36}(t_{36})hw x_1(u_1)\dots x_{36}(u_{36})$,
where $h\in H$, $w$ is a product of some $w_i(1)$.
Then $h_1 x h_1=x_1(\pm t_1)\dots x_{36}(\pm t_{36}) h h' w
x_1(\pm u_1)\dots x_{36}(\pm u_{36})$, and $h_1(-1) x_i(s) h_1(-1)=x_i(s)$
if and only if 
$i=1,2,4,5,6,9,10,13,15,$ $20,26$, $29, 33, 34$, $35,36.$ 

Besides,  $h_1(-1)$ has to commute with $w$. 
As we remember from the previous section,
 $w_{i_1}(1)\dots w_{i_m}(1)$ commutes with $h_1(-1)$,
if and only if $w_{\alpha_{i_1}} \dots w_{\alpha_{i_m}}(\alpha_1)=
\pm \alpha_1$. 

It is clear that any product of $w_1(1)$, $w_2(1)$, $w_4(1)$
$w_5(1)$, $w_6(1)$ and $w_{26}(1)$ commutes with  $h_1(-1)$. 
It is easy to show that other elements of the Weil group
do not commute with $h_1(-1)$. Therefore, the centralizer of 
$h_1(-1)$ is a product of the group $H$ and the subgroup of~$G$, generated
by $X_1$, $X_2$, $X_4$, $X_5$, $X_6$ and $X_{26}$.
The group $X_1$ commutes with all other 
generating subgroups, and the subgroups $X_2$, $X_4$, $X_5$, $X_6$ and $X_{26}$
generate the Chevalley group of the type $A_5$. Thus, the centralizer
of $h_1$ is a product of  $H$ and the central product of
the Chevalley groups $A_1$ and $A_5$, and the commutant of this centralizer
is just a central product of  Chevalley groups $A_1$ and
$A_5$. Consequently, the central quotient of this commutant
is $PSL_2(K)\times PSL_6(K)$.
                      
The centralizer of $h_1h_2$ is generated by $H$ and the subgroup
generated by $X_1$, $X_2$, $X_5$, $X_6$ and $X_8$ (it has the type
 $D_5$), and the central quotient of the commutant of this centralizer
is isomorphic to $PEO_{10}(K)$.

Now we need to separate  Chevalley groups of the type $E_6$
from all classical Chevalley groups. 
It is easy, if we mention that there is no classical Chevalley group
with a central quotient of a commutant of a centralizer of some
its involution
elementary equivalent to
 $PSL_2(K)\times PSL_6(K)$,
and with a central quotient of a commutant of a centralizer of some other
involution elementary equivalent to $PEO_{10}(K)$. 
It can not be true also for the groups of the types $G_2$ and $F_4$.
\end{proof}

\section{Chevalley groups of the type $E_7$.}\leavevmode

In the root system $E_7$ there are $126$ roots, and $63$ of them are
positive.
These positive roots are generated by seven simple roots
$\alpha_1=\frac{1}{2}(e_1+e_8-e_2-e_3-e_4-e_5-e_6-e_7)$, $\alpha_2=e_1+e_2$, 
$\alpha_3=e_2-e_1$, $\alpha_4=e_3-e_2$, $\alpha_5=e_4-e_3$,
$\alpha_6=e_5-e_4$ и $\alpha_7=e_6-e_5$.  
These are other positive roots:
{\small $\alpha_8=1+3$, $\alpha_9=3+4$, $\alpha_{10}=
4+5$, $\alpha_{11}=5+6$, $\alpha_{12}=6+7$,
$\alpha_{13}=1+3+4$,  $\alpha_{14}=2+4$, $\alpha_{15}=3+4+5$,
$\alpha_{16}=2+3+4$,
$\alpha_{17}=2+4+5$, $\alpha_{18}=4+5+6$, $\alpha_{19}=5+6+7$,
$\alpha_{20}=4+5+6+7$, 
$\alpha_{21}=1+2+3+4$, 
$\alpha_{22}=2+3+4+5$, 
$\alpha_{23}=3+4+5+6+7$, 
$\alpha_{24}=2+4+5+6+7$,
$\alpha_{25}=2+3+4+5+6+7$,
$\alpha_{26}=1+3+4+5$,
$\alpha_{27}=3+4+5+6$,
$\alpha_{28}=2+4+5+6$,
$\alpha_{29}=2+3+4+4+5$,
$\alpha_{30}=2+3+4+5+6$,
$\alpha_{31}=1+2+3+4+5$,
$\alpha_{32}=1+3+4+5+6$,
$\alpha_{33}=1+2+3+4+4+5$,
$\alpha_{34}=1+2+3+3+4+4+5$,
$\alpha_{35}=2+3+4+4+5+6$,
$\alpha_{36}=2+3+4+4+5+5+6$,
$\alpha_{37}=1+2+3+3+4+4+5+6$,
$\alpha_{38}=1+2+3+4+5+6$,
$\alpha_{39}=1+2+3+4+4+5+6$,
$\alpha_{40}=1+2+3+4+4+5+5+6$,
$\alpha_{41}=1+2+3+3+4+4+5+5+6$,
$\alpha_{42}=1+2+3+3+4+4+4+5+5+6$,
$\alpha_{43}=1+2+2+3+3+4+4+4+5+5+6$,
$\alpha_{44}=1+1+2+2+3+3+3+4+4+4+4+5+5+5+6+6+7$,
$\alpha_{45}=1+3+4+5+6+7$,
$\alpha_{46}=2+3+4+4+5+6+7$,
$\alpha_{47}=2+3+4+4+5+5+6+7$,
$\alpha_{48}=1+2+3+3+4+4+5+6+7$,
$\alpha_{49}=1+2+3+4+5+6+7$,
$\alpha_{50}=1+2+3+4+4+5+6+7$,
$\alpha_{51}=1+2+3+4+4+5+5+6+7$,
$\alpha_{52}=1+2+3+3+4+4+5+5+6+7$,
$\alpha_{53}=1+2+3+3+4+4+4+5+5+6+7$,
$\alpha_{54}=1+2+2+3+3+4+4+4+5+5+6+7$,
$\alpha_{55}=1+2+3+4+4+5+5+6+6+7$,
$\alpha_{56}=1+2+3+3+4+4+5+5+6+6+7$,
$\alpha_{57}=1+2+3+3+4+4+4+5+5+6+6+7$,
$\alpha_{58}=1+2+2+3+3+4+4+4+5+5+6+6+7$,
$\alpha_{59}=2+3+4+4+5+5+6+6+7$,
$\alpha_{60}=1+2+3+3+4+4+4+5+5+5+6+6+7$,
$\alpha_{61}=1+2+2+3+3+4+4+4+5+5+5+6+6+7$,
$\alpha_{62}=1+2+2+3+3+4+4+4+4+5+5+5+6+6+7$,
$\alpha_{63}=1+2+2+3+3+3+4+4+4+4+5+5+5+6+6+7$.}

\begin{lemma}\label{razlE7}
There exists a first order sentence $\varphi_{E_7}$
that holds in any adjoint groups of the type $E_7$ 
and does not hold in all classical adjoint Chevalley groups,
and also in adjoint Chevalley groups of types $G_2$, $F_2$, $E_6$. 
\end{lemma}

\begin{proof}
Let $Z(h_1)\ni x=
x_1(t_1)\dots x_{63}(t_{63})hw x_1(u_1)\dots x_{63}(u_{63})$,
where $h\in H$, $w$ is a product of some elements $w_i(1)$.
Then $h_1 x h_1=x_1(\pm t_1)\dots x_{63}(\pm t_{63}) h h' w$ $
x_1(\pm u_1)\dots x_{63}(\pm u_{63})$, and $h_1(-1) x_i(s) h_1(-1)=x_i(s)$
if and only if 
$i=1,2,4,5,6,7,10,11,12,14,17,18$, $19,20,24,28,34,37$,
$41,42,43,48,52,53,54,56,57,58,60,61,62$.
Besides, $h_1(-1)$ has to commute with $w$. 
For $w_{i_1}(1)\dots w_{i_m}(1)$ commuting with $h_1(-1)$,
we need $w_{\alpha_{i_1}} \dots w_{\alpha_{i_m}}(\alpha_1)=
\pm \alpha_1$. 

It is clear that any product of $w_1(1)$, $w_2(1)$, $w_4(1)$
$w_5(1)$, $w_6(1)$, $w_7(1)$ and $w_{34}(1)$ commutes with $h_1(-1)$,
and other elements of the Weil group 
does not commute with $h_1(-1)$. Therefore the centralizer
of $h_1(-1)$ is a product of the group $H$ and the subgroup of~$G$, 
generated by $X_1$, $X_2$, $X_4$, $X_5$, $X_6$, $X_7$ and $X_{34}$, consequently,
the central quotient of the commutant of the centralizer
of $h_1(-1)$
is isomorphic to $PSL_2(K)\times PEO_{12}(K)$. Let us look, if
there exists any Chevalley group with a central quotient of
a commutant of a centralizer of some involution is also
isomorphic to this group.

It is clear that it is impossible neither for groups $G_2$, $F_4$, $E_6$,
nor for groups $A_l, C_l$, $D_l$.
The unique possible variant is the group $EO_{15}(K)$ (the group of the 
type~$B_7$).
But in this group there exists an involution
$$
diag[\underbrace{-1,\dots,-1}_{8}, \underbrace{1,\dots,1}_7],
$$
with the central quotient of the commutant of its centralizer isomorphic to
 $PEO_8(K)\times EO_7(K)$.
It can not be possible for a Chevaley group of the type~$E_7$,
since it does not contain a subgroup $EO_7(K)$
(since no roots of the system~$E_7$ can form the subsystem~$B_3$).

Therefore the Chevalley group  $E_7$ 
can not be elementarily equivalent
to the Chevalley group $B_7$, and therefore to any considered
group.
\end{proof}

\section{Chevalley groups of the type $E_8$.}\leavevmode

In the root system $E_8$ there are $240$ roots,  $120$ roots are positive.
These positive roots are generated by eight simple rots
{\small $\alpha_1=\frac{1}{2}(e_1+e_8-e_2-e_3-e_4-e_5-e_6-e_7)$, $\alpha_2=e_1+e_2$, 
$\alpha_3=e_2-e_1$, $\alpha_4=e_3-e_2$, $\alpha_5=e_4-e_3$,
$\alpha_6=e_5-e_4$, $\alpha_7=e_6-e_5$ и $\alpha_8=e_7-e_6$}.  
These are other positive rots (for short we omit the signs ``$+$''):
{\footnotesize 
$\alpha_9=13$, $\alpha_{10}=24$, $\alpha_{11}=
34$, $\alpha_{12}=45$, $\alpha_{13}=56$,
$\alpha_{14}=67$, $\alpha_{15}=78$, $\alpha_{16}=134$,  $\alpha_{17}=234$, 
$\alpha_{18}=245$,
$\alpha_{19}=1234$,
$\alpha_{20}=345$, $\alpha_{21}=456$, $\alpha_{22}=1345$,
$\alpha_{23}=567$, 
$\alpha_{24}=2456$, 
$\alpha_{25}=3456$, 
$\alpha_{26}=678$, 
$\alpha_{27}=4567$,
$\alpha_{28}=13456$,
$\alpha_{29}=2345$,
$\alpha_{30}=23445$,
$\alpha_{31}=23456$,
$\alpha_{32}=12345$,
$\alpha_{33}=23567$,
$\alpha_{34}=123445$,
$\alpha_{35}=123456$,
$\alpha_{36}=34567$,
$\alpha_{37}=1233445$,
$\alpha_{38}=45678$,
$\alpha_{39}=234456$,
$\alpha_{40}=1234456$,
$\alpha_{41}=134567$,
$\alpha_{42}=5678$,
$\alpha_{43}=234567$,
$\alpha_{44}=1234567$,
$\alpha_{45}=245678$,
$\alpha_{46}=2344556$,
$\alpha_{47}=12344556$,
$\alpha_{48}=345678$,
$\alpha_{49}=12334456$,
$\alpha_{50}=1345678$,
$\alpha_{51}=2344567$,
$\alpha_{52}=12344567$,
$\alpha_{53}=2345678$,
$\alpha_{54}=23445567$,
$\alpha_{55}=123344556$,
$\alpha_{56}=23445566778$,
$\alpha_{57}=23445678$,
$\alpha_{58}=1233444556$,
$\alpha_{59}=234455667$,
$\alpha_{60}=12345678$,
$\alpha_{61}=123445567$,
$\alpha_{62}=1234455667$,
$\alpha_{63}=123445678$,
$\alpha_{64}=12233444556$,
$\alpha_{65}=123344567$,
$\alpha_{66}=1233445567$,
$\alpha_{67}=12334445567$,
$\alpha_{68}=12334455667$,
$\alpha_{69}=1233445678$,
$\alpha_{70}=122334445567$,
$\alpha_{71}=234455678$,
$\alpha_{72}=2344556678$,
$\alpha_{73}=1234455678$,
$\alpha_{74}=12344556678$,
$\alpha_{75}=123344455678$,
$\alpha_{76}=123344455667$,
$\alpha_{77}=12334455678$,
$\alpha_{78}=123445566778$,
$\alpha_{79}=123344556678$,
$\alpha_{80}=1233445566778$,
$\alpha_{81}=1233444556678$,
$\alpha_{82}=12334445566778$,
$\alpha_{83}=1233444555667$,
$\alpha_{84}=12334445556678$,
$\alpha_{85}=123344455566778$,
$\alpha_{86}=1233444555666778$,
$\alpha_{87}=1223344455678$,
$\alpha_{88}=1223344455667$,
$\alpha_{89}=12233444556678$,
$\alpha_{90}=122334445566778$,
$\alpha_{91}=12233444555667$,
$\alpha_{92}=1223344455566778$,
$\alpha_{93}=1223344455566778$,
$\alpha_{94}=12233444555666778$,
$\alpha_{95}=122334444555667$,
$\alpha_{96}=122334444555666778$,
$\alpha_{97}=12233444455566778$,
$\alpha_{98}=122334444555666778$,
$\alpha_{99}=1223344445555666778$,
$\alpha_{100}=1223334444555667$,
$\alpha_{101}=12233344445556678$,
$\alpha_{102}=122333444455566778$,
$\alpha_{103}=1223334444555666778$,
$\alpha_{104}=12233344445555666778$,
$\alpha_{105}=122333444445555666778$,
$\alpha_{106}=1222333444445555666778$,
$\alpha_{107}=11223334444555667$,
$\alpha_{108}=112233344445556678$,
$\alpha_{109}=1122333444455566778$,
$\alpha_{110}=11223334444555666778$,
$\alpha_{111}=112233344445555666778$,
$\alpha_{112}=1122333444445555666778$,
$\alpha_{113}=11223333444445555666778$,
$\alpha_{114}=11222333444445555666778$,
$\alpha_{115}=112223333444445555666778$,
$\alpha_{116}=1122233334444445555666778$,
$\alpha_{117}=11222333344444455555666778$,
$\alpha_{118}=112223333444444555556666778$,
$\alpha_{119}=1122233334444445555566667778$,
$\alpha_{120}=11222333344444455555666677788$.}

\begin{lemma}\label{razlE8}
There exists a first order sentence $\varphi_{E_8}$
that holds in any adjoint group of the type $E_8$ and does not
hold in all adjoint Chevalley groups of other types.
\end{lemma}

\begin{proof}
Consider the centralizer of $h_1$.

Let $Z(h_1)\ni x =x_1(t_1)\dots x_{120}(t_{63})h
w x_1(u_1)\dots x_{120}(u_{120})$,
where $h\in H$, $w$ is a product of some  $w_i(1)$.
Then $h_1 x h_1=x_1(\pm t_1)\dots x_{120}(\pm t_{120}) h h' w$
$x_1(\pm u_1)\dots x_{120}(\pm u_{120})$, and $h_1(-1) x_i(s) h_1(-1)=x_i(s)$
if and only if 
$i=1,2,4,5,6,7,8,10,11,12,$ $13,14,15,18,19,21,23,24,$
$26, 27,33,37,38,$
$42,45,49,55,58,64,65,66$, $67,68$, $69,70$, $75,76,77,79,80,$
$81,82,83,84$, $85,86,87$, $88,89,90$,
$91,92,93,94,95,96,$
$97,98$, $99,113$, $115,116,$
$117,118,119,120.$

It is clear that any product of $w_1(1)$, $w_2(1)$, $w_4(1)$
$w_5(1)$, $w_6(1)$, $w_7(1)$, $w_8(1)$ and $w_{37}(1)$ commutes with $h_1(-1)$, 
and other elements of the Weil group do not commute, therefore
the centralizer of $h_1(-1)$ is a product of the group $H$ and the subgroup
of~$G$, generated by $X_1$, $X_2$, $X_4$, $X_5$, $X_6$, $X_7$, $X_8$ and
 $X_{37}$.

The central quotient of the commutant of the centralizer of this involution
is isomorphic to $PSL_2(K)\times G(E_7,K)$, therefore the group of the type
 $E_8$ can not be elementary equivalent to any other Chevalley group.
\end{proof}

Lemmas \ref{klass1}--\ref{razlE8} imply

\begin{proposition}\label{root}
If two \emph{(}elementary\emph{)} Chevalley groups over infinite fields
of characteristic $\ne 2$ are elementarily equivalent, then the corresponding
root systems coincide.
\end{proposition}

\section{Definability of a field in Chevalley groups.}\leavevmode

To proof that in every Chevalley group its field is definable
we only need to proof that in any adjoint Chevalley group of the type
 $A_1$ (i.e., in $PSL_2(K)$) its basic field is definable.

\begin{lemma}\label{fieldA1}
If the groups $PSL_2(K)$ and $PSL_2(K')$ \emph{(}$K,K'$ are infinite fields
of characteristics $\ne 2$\emph{)} are elementarily equivalent, then 
the fields $K,K'$ are elementarily equivalent.
\end{lemma}
\begin{proof}
Suppose that we have the group $PSL_2(K)$.

Consider matrices satisfying the formula
$$
Cell(M):= M^2\ne 1\land 
\exists X ((XMX^{-1})M=M(XMX^{-1})\land X^2M\ne MX^2).
$$
Show that these matrices in some basis have the form
$$
A=\pm \begin{pmatrix}
1& t\\
0& 1
\end{pmatrix}.
$$

Let $M$ be diagonalizable in the field $\overline K$, then in 
some basis it has the form
$$
A=\begin{pmatrix}
\alpha& 0\\
0& 1/\alpha
\end{pmatrix}.
$$
A matrix of such a form can commute either with diagonal matrix, or
we can have the situation
$AB=-BA$ in the group $SL_2(K)$, i.e. 
$$
A=\begin{pmatrix}
i& 0\\
0& -i
\end{pmatrix},
$$
But in this case $A^2=1$ in the group $PSL_2(K)$, therefore,
$M^2=1$, but it contradicts to the first condition of the formula $Cell$.

Any matrix $M$, conjugate to the matrix of the form
$$
A=\pm \begin{pmatrix}
1& t\\
0& 1
\end{pmatrix},
$$
satisfies the formula $Cell$ (we can choose $X$ diagonal with distinct
elements on the diagonal, but with squares  $\ne \pm 1$).

Fix some matrix $A$, satisfying the formula $Cell$. Consider
all matrices commuting with~$A$.
Let for a given~$A$ this set be denoted by ${\cal X}_A$.
These are matrices, that in the same basis have the form
$$
\begin{pmatrix}
1& t\\
0& 1
\end{pmatrix},\quad t\in K.
$$

The normalizer $N({\cal X}_A)$ of ${\cal X}_A$ consists of matrices
having in the given basis the form
$$
\begin{pmatrix}
\alpha& \beta\\
0& 1/\alpha
\end{pmatrix}.
$$

Let us choose one matrix  $B$ from $N({\cal X}_A)\setminus {\cal X}_A$
and all matrices commuting with it.
Denote this set by ${\cal D}_{A,B}$.
We can find a basis, where all matrices from ${\cal D}_{A,B}$
are
$$
\begin{pmatrix}
\alpha& 0\\
0& 1/\alpha
\end{pmatrix},
$$
and the matrices from ${\cal X}_A$ are
$$
\pm \begin{pmatrix}
1& t\\
0& 1
\end{pmatrix}.
$$

Besides it, consider the set ${\cal E}_{A,B}$ of matrices, normalizing
${\cal D}_{A,B}$, but not belonging to ${\cal D}_{A,B}$.
They are
$$
\begin{pmatrix}
0& \gamma\\
-1/\gamma& 0
\end{pmatrix}.
$$

By  ${\cal Y}_{A,B}$ we denote the matrices conjugate to
the matrices from ${\cal X}_A$ by matrices from ${\cal E}_{A,B}$.
They are
$$
\pm \begin{pmatrix}
1& 0\\
s& 1
\end{pmatrix}.
$$

Now let us introduce the formula
$$
\varphi(M,N):= (M\in {\cal X}_A) \land (N\in {\cal Y}_{A,B}) \land
(MNM\in {\cal E}_{A,B}).
$$

If for a matrix
$$
M=\begin{pmatrix}
1& t\\
0& 1
\end{pmatrix}
$$
we have $\varphi(M,N)$ for some~$N$, then
$$
N=M^\varphi=\begin{pmatrix}
1& 0\\
-1/t& 1
\end{pmatrix}.
$$

Let us write 
$$
E_t:= {\cal E}(A_t)=A_tA_t^\varphi A_t=
\begin{pmatrix}
0& t\\
-1/t& 0
\end{pmatrix}.
$$

Fix 
$$
A_u=\begin{pmatrix}
1& u\\
0& 1
\end{pmatrix}.
$$

Since
$$
E_tE_u^{-1}=
\begin{pmatrix}
0& t\\
-1/t& 0
\end{pmatrix}
\begin{pmatrix}
0& -u\\
1/u& 0
\end{pmatrix}=
\begin{pmatrix}
t/u& 0\\
0& u/t
\end{pmatrix},
$$
we have 
$$                     
E_{ts/u}=E_t E_u^{-1} E_s.
$$

Note that since we have the group $PSL_2(K)$, the mapping ${\cal 
E}: {\cal X}_A\to {\cal E}_{A,B}$, ${\cal X}_A\mapsto {\cal E}({\cal X}_A)$, 
is not invertible (namely, it is not monomorphic), 
because ${\cal E}(A_t)={\cal E}(A_{-t})$
in the group $PSL_2(K)$.

Suppose that we have $A_t$ and $A_s$ and want to find
 $A_t\otimes A_s=A_{ts/u}$.
In this case the mapping $A_t\mapsto t/u$ is an isomorphism between
the set ${\cal X}_A$ (with the operation of addition
 $A_t\oplus A_s=A_tA_s$ and the operation
of multiplication~$\otimes$) and the field~$K$. Since
 $E_{ts/u}={\cal E}(A_t){\cal E}(A_u)^{-1}
{\cal E}(A_s)$, we only need to restore $A_{ts/u}$ by
$E_{ts/u}$. We have only to candidates:
 $A_{ts/u}$ and $A_{-ts/u}$, and we need to choose the right one.

There are the following conditions:
$$
A_tA_{ts/u}=A_{t(1+s/u)}=A_{t/u(s+u)}=A_t\otimes (A_sA_u),
$$
therefore,
$$
{\cal E}(A_tA_{ts/u})={\cal E}(A_t){\cal E}(A_u)^{-1} {\cal E}(A_sA_u)
$$
for $t,s\ne 0$.

Also we have
$$
A_tA_{-ts/u}=A_{t/u(u-s)}.
$$
The condition 
$$
{\cal E}(A_tA_{-ts/u})={\cal E}(A_t){\cal E}(A_u)^{-1} {\cal E}(A_sA_u)
$$
for $t,s\ne 0$, as it is easy to check, can not hold, therefore
the formula of multiplication can be written
in the next form:
\begin{multline*}
A_w=A_t\otimes A_s \Longleftrightarrow (A_t=1\land A_w=1)\lor (A_s=1\land A_w=1)\lor\\
\lor (A_s\ne 1 \land A_t\ne 1 \land {\cal E}(A_w)={\cal E}(A_t){\cal E}(A_u)^{-1}{
\cal E}(A_s)\land {\cal E}(A_tA_w)={\cal E}(A_t){\cal E}(A_u)^{-1} {\cal E}(A_sA_u)).
\end{multline*}

Consequently, the isomorphism is constructed, thus,
the field $K$ is definable
in the group $PSL_2(K)$, i.\,e. there exist formulas
 $\varphi(X_1,X_2,\dots,X_l;Y)$,
$\psi(X_1,\dots,X_l)$, $Add(Y_1,Y_2,Y_3)$, $Mult(Y_1,Y_2,Y_3)$ such that
or any $X_1,\dots,X_l$, satisfying the formula $\psi(\dots)$,
all such $Y$ that  $\varphi(X_1,\dots,X_l;Y)$ holds, with the operations
of addition and multiplication, given by the formulas
 $Add(\dots)$ and $Mult(\dots)$,
form a field, isomorphic to~$K$.

Therefore, if the groups $PSL_2(K)$ and $PSL_2(K')$ are elementarily 
equivalent, then the fields $K,K'$ are elementarily equivalent.
\end{proof}

\begin{proposition}\label{fields}
If two adjoint elementary Chevalley groups $G(\Phi,K)$
and $G(\Phi,K')$  \emph{(}$K,K'$ are infinite fields of characteristics 
$\ne 2$\emph{)}
are elementarily equivalent, then the fields $K,K'$ are elementarily 
equivalent.
\end{proposition}
\begin{proof}
For a given root system~$\Phi$ consider the central quotient of the
commutant of the centralizer of a suitable involution
 (see Lemma~\ref{klass1}--\ref{razlE8}) and
choose there such a component of the direct product that
is isomorphic to $PSL_2(K)$. Then we use Lemma~\ref{fieldA1}.
\end{proof}

\section{Isomorphism of weight lattices.}\leavevmode
 
Now we only need to prove

\begin{proposition}\label{lattice}
If two \emph{(}elementary\emph{)} Chevalley 
groups are elementarily equivalent, then their weight lattices
coincide.
\end{proposition}
\begin{proof}
We shall consider elementary Chevalley groups.

Note that if two Chevalley groups are elementarily equivalent,
then their centers are isomorphic, since they are finite.
The center of a Chevalley group is a quotient group of its weight lattice 
by the root lattice. We know that two elementary
Chevalley group have the same type. Then the quotient group of the weight
lattice by the root lattice completely defines the weight lattice,
except the single case:
the Chevalley group of the type $D_{2m}$, $m\ge 3$, and
intermediate lattice such that the order of the quotient 
of weight lattice by the root lattice is~$2$. In this case there are
two non-isomorphic Chevalley groups: one of them is $SO_{4m}(K)$ 
(commutant of this group in elementary case), 
and the other is a so-called ``semi-spinor group''.

Recall (see~\cite{Burbaki}), that roots of the system 
$D_{2m}$ are the vectors $\pm e_i\pm e_j$, $i\ne j$ 
in the othonormal basis $e_1,\dots, e_{2m}$ of the $2m$-dimensional
Euclidean space, the weight lattice of this system consists of
the vectors with even sum of integral coordinates,
the universal lattice is
$$
\oplus_{i=1}^{2m} \mathbb Z e_i + \mathbb Z\left( \frac{1}{2} \sum_{i=1}^{2m} e_i\right),
$$
the first (``orthogonal'') intermediate lattice $\Lambda_1$ is
$$
\oplus_{i=1}^{2m} \mathbb Z e_i,
$$
the second (``semi-spinor'') $\Lambda_2$ is
$$
\oplus_{i=1}^{2m-1} \mathbb Z (e_i-e_{i+1}) \oplus \mathbb Z (e_{2m-1}+e_{2m})+
\mathbb Z\left( \frac{1}{2} \sum_{i=1}^{2m} e_i\right).
$$

Suppose that we have the group of the type
$D_{2m}$, $m\ge 4$. Then we can consider (define by a formula) 
an involution~$j$ with a central quotient of the commutant
of the centralizer having the type $D_6\times D_{2m-6}$ 
(for example, it is the involution $h_{\alpha_1}(-1)
h_{\alpha_3}(-1)h_{\alpha_5}(-1)$). 
The commutant of the centralizer of~$j$ is a central product
of the group~$D_6$ and the group~$D_{2m-6}$, and the weight lattice
of the subgroup~$D_6$ is the same as the initial group has 
 (if the initial group is orthogonal, then the subgroup 
is orthogonal, if the initial group is semi-spinor, then  
the subgroup is semi-spinor). Therefore we shall take the subgroup of the 
type~$D_6$ to define its lattice type.

Let us have the group~$G$ of the type $D_6$  with intermediate lattice.
We shall suppose that either $G=G_1=SO_{12}(K)'$, or $G=G_2$ is a semi-spinor
group (that is the quotient group of $Spin_{12}(K)$ by the central
subgroup having two elements).

The center of each of the groups $G_1$, $G_2$ has two elements and 
the central quotient is the adjoint group $PSO_{12}(K)'$, 
that is already studied.

The subgroup $H_1=H(G_1)$ is generated by
$h_{\alpha_1}(t_1)$, $h_{\alpha_2}(t_2)$,
$h_{\alpha_3}(t_3)$, $h_{\alpha_4}(t_4)$, $h_{\alpha_5}(t_5)$,
$h_{\alpha_6}(t_6)$ 
with the condition $h_{\alpha_5}(-1)h_{\alpha_6}(-1)=1$,
and the subgroup $H_2=H(G_2)$ is generated by the same elements
with the condition
$$
h_{\alpha_1}(-1)h_{\alpha_3}(-1)h_{\alpha_5}(-1)=1.
$$

In the group $PSO_{12}(K)$ for $i\notin K$ in the subgroup~$H$ there are
 16 involutions, 
that are generated by $h_{\alpha_1}(-1), h_{\alpha_2}(-1), h_{
\alpha_3}(-1), h_{\alpha_4}(-1)$. The maximal set of commuting 
involutions consists of $2^9$ involutions, 
that are generated by the mentioned elements, and also by elements
$w_{e_1-e_2}w_{e_1+e_2}$, $w_{e_2-e_3}w_{e_2+e_3}$, 
$w_{e_3-e_4}w_{e_3+e_4}$, $w_{e_4-e_5}w_{e_4+e_5}$, 
$w_{e_5-e_6}w_{e_5+e_6}=w_{\alpha_5}w_{\alpha_6}$. 
If  $i\in K$, then we can also add the involution
$h_{\alpha_5}(i)h_{\alpha_6}(i)$, and therefore we have the set of 
$2^{10}$ (there can not be more involutions in the group $PSO_{12}(K)$!) 
involutions.

Now we can note that in the group $G_1=SO_{12}(K)$ the inverse images of these 
involutions also are involutions, since in this group we have the condition
 $h_{e_i-e_{i+1}}(-1)=h_{e_i+e_{i+1}}(-1)$, 
$i=1,\dots,5$. Besides that, we can add the central involution
$-1=h_{\alpha_1}(-1)h_{\alpha_3}(-1)h_{\alpha_5}(-1)$. 

In the semi-spinor group $G_2$ the condition
$h_{e_i-e_{i+1}}(-1)=h_{e_i+e_{i+1}}(-1)$ does not hold, 
therefore $(w_{e_i-e_{i+1}}w_{e_i+e_{i+1}})^2=h_{e_i-e_{i+1}}(-1) 
h_{e_i+e_{i+1}}(-1)\ne 1$. Similarly, 
$(h_{\alpha_5}(i)h_{\alpha_6}(i))^2=h_{\alpha_5} (-1) 
h_{\alpha_6}(-1)\ne 1$. Consequently, the inverse images of only
16 involutions from the adjoint group are involutions 
in the semi-spinor group.

Thus, we can easily write a sentence, separating the groups
$G_1$ and $G_2$, i.\,e., they can not be elementarily equivalent.

Therefore, elementary equivalent Chevalley groups has the same 
weight lattices.
\end{proof}

Theorem \ref{easytheor} and Propositions \ref{elemsubgr}, \ref{root},
\ref{fields},~\ref{lattice} imply the main

\begin{theorem}
Let 
 $G=G_\pi (\Phi,K)$ и $G'=G_{\pi'}(\Phi',K')$
\emph{(}or $E_\pi (\Phi,K)$ and $E_{\pi'}(\Phi',K'))$ be two
\emph{(}elementary\emph{)}
Chevalley groups over infinite fields $K$ and~$K'$ of characteristics~$\ne 2$,
with weight lattices $\Lambda$ and $\Lambda'$, respectively.
Then the groups $G$ and $G'$ are elementarily equivalent 
if and only if the root systems $\Phi$ and $\Phi'$ are isomorphic, 
the fields $K$  and $K'$ 
are elementarily equivalent, the lattices $\Lambda$ and $\Lambda'$ coincide.
\end{theorem}


\begin{thebibliography}{99}
\bibitem{KeisChang} 
Chang~C.C. and Keisler~H.J., Model Theory,
North--Holland, Amsterdam--London, American Elsevier, New York, 1973.

\bibitem{Maltsev} 
Maltsev~A.I., On elementary properties of linear groups,
in: Problems of Mathematics and Mechanics [in Russian],
Novosibirsk, 1961, 110--132.


\bibitem{BeiMikh}
 Beidar C.I. and Michalev A.V. On Malcev's theorem on elementary 
equivalence of linear groups. 
Contemporary mathematics, 1992, 131, 29--35.

\bibitem{Bun1}
Bunina~E.I., Elementary equivalence of unitary linear
groups over rings and skewfields, \textit{Russian Mathematical Surveys},
\textbf{53}, No.~2, 137--138 (1998).

\bibitem{Bun2}
Bunina~E.I., Elementary equivalence of Chevalley groups,
\textit{Russian Mathematical Surveys}, \textbf{56}, No.~1, 157--158 (2001).

\bibitem{Bun3} Bunina E.I., Mikhalev A.V. 
Combinatorial and Logical Aspects of Linear Groups and Chevalley Groups. 
Acta Applicandae Mathematicae, 2005, 85(1-3), 57-74.

\bibitem{BunChev}
Bunina~E.I., Chevalley groups over fields and their elementary 
properties,
\textit{Russian Mathematical Surveys}, \textbf{59}, No.~5, 952--953 (2004).

\bibitem{Hamfris} Humphreys J.E. Introduction to Lie algebras and representation
theory. Springer--Verlag New York, 1978.

\bibitem{Burbaki} Bourbaki N. Groupes et Alg\'ebres de Lie. Hermann, 1968.

\bibitem{Steinberg} 
Steinberg R. Lectures on Chevalley groups, Yale University, 1967.

\bibitem{Chevalley} Chevalley C. Certain schemas des groupes semi-simples. 
Sem. Bourbaki,  1960--1961, 219, 1--16.

\bibitem{Artem_dis} Golubkov A.Ju. The prime radical of classical groups over associative rings, ph.d. thesis.
Moscow, 2001.

\bibitem{VavPlotk1}  Vavilov N.A., Plotkin E.B. Chevalley groups over commutative rings. I.
Elementary calculations. Acta Appl. Math. 45 (1996), 73--113.
	
\bibitem{Abe1} Abe E., \emph{Chevalley groups over local rings}, Tohoku Math. J. 21 
(1969), N3, 474--494

\bibitem{A3} Abe E., \emph{Whitehead groups of Chevalley groups over Laurent 
polynomial rings}, Comm. algebra 11 (1983), N12, 1271--1308.

\bibitem{A4} Abe E., \emph{Whitehead groups of Chevalley groups over Laurent 
polynomial rings}, Preprint Univ. Tsukuba (1988).


\bibitem{A6} Abe E., \emph{Normal subgroups of Chevalley groups over 
commutative rings}, Contemp. Math. 83 (1989), 1--17.

\bibitem{A7} Abe E., \emph{Automorphisms of Chevalley groups over 
commutative rings}, Алгебра и анализ, 5 (1993), N3, 74--90.

\bibitem{AH} Abe E., Hurley J., \emph{Centers of Chevalley groups over commutative rings},
Comm. Algebra 16 (1988), N1, 57--74.

\bibitem{AS} Abe E., Suzuki K., \emph{On normal subgroups of Chevalley groups 
over commutative rings}, Tohoku Math. J. 28 (1976), N1, 185--198

\bibitem{BMS} Bass H., Milnor J., Serr J.-P., \emph{Solution of the congruence subgroups
problem for $SL_n$ ($n\ge 3$) and $Sp_{2n}$ ($n\ge 2$)}, Publ. Math. Inst. Hautes Et. Sci.
(1967), N33, 59--137.

\bibitem{artem17} Borel A., \emph{Properties and linear representations of Chevalley
groups}, Lect. Notes in Math., 131, Berlin; Heidelberg; N.Y.: Springer, 1970, 1--55.

\bibitem{Cn} Cohn P., \emph{On the structure of the $GL_2$ of a ring},
Publ. Math. Inst. Hautes Et. Sci. (1966), N30, 365--413.

\bibitem{Cs} Costa D.L., \emph{Zero-dimensionality and the $GE_2$ 
of polynomial rings}, J. Pure Appl. Algebra
50 (1988), 223--229.


\bibitem{Ko1} Kopeiko V.I., \emph{The stabilization of symplectic 
groups over a polynomial ring}, Math. U.S.S.R. Sbornik 34 (1978), 655-669.

\bibitem{Ko2} Kopeiko V.I. \emph{Unitary and orthogonal groups over rings
with involution}, Algebra and discrete Math., Калмытский гос. университет, Элиста,
1985, 3--14.

\bibitem{M} Matsumoto H., \emph{Sur les sous-groupes arithm\'etiques des groupes
semi-simples deploy\'es}, Ann. Sci. Ecole Norm. Sup., 4eme ser. 2 (1969), 1--62.

\bibitem{GMV} Gr\"unewald F., Mennicke J., Vaserstein L.N., \emph{On symplectic 
groups over polynomial rings}, Math. Z. 206 (1991), N1, 35--56.


\bibitem{Su1} Suslin A.A., \emph{On a theorem of Cohn}, J. Sov. Math. 17 (1981), N2,
1801--1803.

\bibitem{Su2} Suslin A.A., \emph{On the structure of geneal linear group over polynomial
ring}, Soviet Math. Izv. 41 (1977), N2. 503--516.

\bibitem{SK} Suslin A.A., Kopeiko V.I., 
\emph{Quadratic modules and orthogonal groups over
polynomial rings}, J.Sov. Math. 20 (1982), N6, 2665--2691.

\bibitem{Sw} Swan R., \emph{Generators and relations for certain special linear groups},
Adv. Math. 6 (1971), 1--77.

\bibitem{Vo} Vorst T., \emph{The general linear group of polynomial rings over regular 
rings}, Comm. Algebra 9 (1981), N5, 499--509.

\bibitem{Post} Postnikov M.M. Lectures in geometry. Lie groups and Lie algebras.
Semester V, Nauka, 1982 (in Russian).

\bibitem{St3} Stein M.R. \emph{Surjective stability in dimension $0$ for $K_2$
and related functors}, Trans. Amer. Soc., 1973, 178(1), 165--191.

\end{thebibliography}
\end{document}